\documentclass[10pt]{article}
\input epsf
\usepackage{enumerate, fullpage}
\usepackage{graphicx}
\usepackage{multirow}
\usepackage{color}
\usepackage{amsmath}
\usepackage{amsfonts}
\usepackage{amssymb}
\usepackage{mathrsfs}
\usepackage{lscape}

\def\cF{{\cal F}}
\def\cB{{\cal B}}

\def\abs#1{|#1|}
\def\norm#1{\|#1\|}

\newcommand{\R}{\mathbb{R}}
\newcommand{\F}{\mathbb{F}}
\newcommand{\Z}{\mathbb{Z}}
\newcommand{\closure}{{\mathop{\rm cl}}}
\newcommand{\interior}{\mathop{\rm int}}
\newcommand{\conv}{\mathop{\rm conv}}

\newcommand{\bmath}[1]{ {\boldmath $ {#1}$}}
\newcommand{\diag}{\mathop{\rm diag}}
\newcommand{\beq}{\begin{equation}}
\newcommand{\eeq}{\end{equation}}

\newcommand{\beqnr}{\begin{eqnarray}}
\newcommand{\eeqnr}{\end{eqnarray}}

\newcommand{\benum}{\begin{enumerate}}
\newcommand{\eenum}{\end{enumerate}}

\newcommand{\argmin}{\mathop{\rm argmin}}
\newcommand{\QED}{\rule{.1in}{.1in}}

\newcommand{\cH}{{\cal H}}

\newcommand{\cP}{{\cal P}}

\newtheorem{DE}{Definition}[section]

\newtheorem{PR}[DE]{Proposition}

\newtheorem{LE}[DE]{Lemma}
\newtheorem{EX}[DE]{Example}
\newtheorem{RE}[DE]{Remark}
\newtheorem{THM}[DE]{Theorem}
\newtheorem{CO}[DE]{Corollary}

\newcommand{\qed}{\mbox{}\hspace*{\fill}\nolinebreak\mbox{$\rule{0.7em}{0.7em}$}}

\advance \topmargin by -\headheight
\advance \topmargin by -\headsep

\evensidemargin \oddsidemargin

\begin{document}

\begin{center}
\Large {\bf Cutting-planes for optimization of convex functions over nonconvex sets}

\vskip 0.1cm

\normalsize

Daniel Bienstock and Alexander Michalka, Columbia University \\
\vskip 0.1cm
November 2011, revised May 2013
\end{center}

\begin{abstract}
We derive linear inequality characterizations for sets
of the form $\conv \{ (x, q) \, \in \R^d \times \R  \, : \, q \ge Q(x), \ \ x \in \R^d - \interior(P) \}$ where $Q$ is convex and differentiable and $P \subset \R^d$.  We show that in several cases our characterization leads
to polynomial-time separation algorithms that operate in the original space
of variables, in particular
when $Q$ is a positive-definite quadratic and $P$ is a 
polyhedron or an ellipsoid.

\end{abstract}

\section{Introduction}
The current state-of-the-art for linear mixed-integer programming relies on
cutting-planes, a methodology
supported by a strong body of theory that has also achieved
computational success. Nevertheless, the solution of an optimization problem 
$\min \{ \, Q(x) \, : \, x \in \cF \, \} $
with $Q(x)$ convex and $\cF \subseteq \R^d$ mixed-integer
would present a challenge to the cutting-plane approach.  Any algorithm
that relies on separation from $\conv(\cF)$
will in general fail, because an optimal solution $x^*$ to
$\min \{ \, Q(x) \, : \, x \in \conv(\cF) \, \}$ may 
satisfy (i) $x^* \notin \cF$, and (ii) $x^*$ is
in the relative
interior of a face of $\conv(\cF)$, and so no cutting-plane can separate 
$x^*$ from $\cF$.  \\

\noindent This observation suggests a paradigm used in the
``lattice-free set'' methodology in mixed-integer programming (reviewed below). Given $(x^*,q^*) \in \R^d\times \R$ with $x^* \notin \cF$, one computes a set $P \subset \R^d$ with $x^* \in \interior(P)$ and $\cF \cap \interior( P) = \emptyset$ (``int'' denotes interior) together with an inequality that separates $(x^*, q^*)$
from the set
\begin{eqnarray}
S & \doteq & \{ \, (x, q) \, \in \R^d \times \R  \ : \ q \ge Q(x), \ \ x \in \R^d - \interior(P) \, \}.  \label{Sdef}
\end{eqnarray}

The focus of this paper is the study of sets of the general form (\ref{Sdef}) with the goal of characterizing $\conv(S)$ by linear inequalities.\footnote{Note that for any extreme point $(x, q)$ of $\conv(S)$ we have $x \in S$ and $q = Q(x)$; a folklore observation.}
   Two classes of `trivial'
valid inequalities for $S$ are (a)
valid inequalities for $\R^d - P$; and  (b) linearization, or first-order, inequalities 
\begin{eqnarray}
q & \ge & Q(y) + \nabla Q(y)^T (x - y), \label{unlifted}
\end{eqnarray}
where $y \in \R^d$.  Usually these two families of inequalities
are not sufficient to characterize $conv(S)$.  In this paper,
motivated by mixed-integer programming considerations,  we consider 
'lifted' versions of (\ref{unlifted}), that is to say inequalities of the form
\begin{eqnarray}
q & \ge & Q(y) + \nabla Q(y)^T (x - y) + \alpha \, p^T  (x - y), \label{firstlifted0}
\end{eqnarray}
where $\alpha > 0$ and $p \in \R^d$.  For (\ref{firstlifted0}) to be valid
$y$ must lie in the boundary of $P$; further $p$ cannot be arbitrary, and instead must
point ``into'' $P$ in a sense made precise later.  We obtain the following results:\\

\noindent {\bf Theorem I.} \hspace{.1in} Let $Q(x)$ be convex and differentiable. 
Any linear inequality  $\delta q \ge \beta^T x + \beta_0$ valid for $S$ 
and such that $\{ (x, q) \in \R^d \times \R \, : \, \delta q = \beta^T x + \beta_0 \}$ is a supporting hyperplane for $\conv(S)$ is dominated
by a combination of up to two inequalities of three types: (a) valid 
inequalities for $\R^d - P$, (b) linearization inequalities obtained 
at points $y \in \R^d - \interior(P)$, and (c) valid lifted inequalities obtained at
points $y$ in the boundary of $P$.  \\

\noindent Separation over the three types inequalities listed in Theorem I is closely related, but not
precisely equivalent to separation from $\conv(S)$.  In this regard, we obtain a sharpening of Theorem I:\\

\noindent {\bf Theorem II.}
\hspace{.1in} Suppose that $Q(x)/\norm{x} \rightarrow +\infty$ as $\norm{x} \rightarrow + \infty$, that there is a polynomial-time
separation oracle for $\R^d - \interior(P)$, and that $\nabla Q(x)$ is polynomial-time computable at any
$x$.  Then polynomial-time separation over $\conv(S)$ is
equivalent to polynomial-time separation over the lifted inequalities.\\

\noindent Unlike what happens in the linear mixed-integer setting, one can produce examples where
a lifted inequality (\ref{firstlifted0}) is binding 
at just one point -- $y$.   We show that such
cases can be essentially characterized in terms of the structure of the boundary of $P$ where lifting is attempted and (again) the degree of 
strong convexity of $Q(x)$:\\

\noindent {\bf Theorem III.}
\hspace{.1in}Let $y$ be a point in the boundary of $P$ such that an open half-ball with center $y$ and positive radius is contained in $P$.  Suppose
further that $Q(x)/\norm{x} \rightarrow +\infty$ as $\norm{x} \rightarrow + \infty$.  Then any lifted inequality (\ref{firstlifted0}) obtained at $y$, and using the maximum valid lifting
coefficient $\alpha$, will be binding at some point $(w, Q(w))$ where $w \neq y$ is in the boundary of $P$. \\

\noindent Finally, in several cases the characterization provided by Theorem I leads to polynomial-time separation: \\ 

\noindent {\bf Theorem IV.}
\hspace{.1in}One can separate in polynomial time from:
\begin{itemize}
\item[(i)] A set $\conv(S)$ as above, when $Q(x)$ is a positive-definite quadratic and
$P$ is a polyhedron or an ellipsoid.
\item[(ii)] A set of the form $\{(x, w, q) \in \R^n \times \R \times \R \, : \, q \ge x^THx + h^Tx, \, w \le x^TAx \}$, where $H \succ 0$ and $A \succeq 0$. 
\end{itemize}

\noindent Case (i) is important because the exclusion of a polyhedron 
has been proposed in several of the lattice-free set schemes in the literature. The ellipsoidal case arises, for example, when considering the cardinality-constrained
convex quadratic programming problem \cite{danoeig}; also see \cite{bucalo}.
In Section \ref{twoquad} we provide motivation for the study of the set in (ii); however, using
$d = n + 1$ and $P = \{(x,w) \in \R^n \times \R \, : \, w \geq x^TAx \}$ 
this set is of the form (\ref{Sdef}).

\subsubsection{Motivation and background}
The ``lattice-free set'' paradigm can
be considered one of the single most fundamental ideas underlying 
the theory of cutting-planes for linear mixed-integer programming.  Our work in this
paper seeks to extend this methodology to the nonlinear setting.
In the linear, pure integer case the methodology can be outlined
as follows.  Let $\cF = \{ x \in \Z^d \, : \, Ax \ge b \, \}$ and consider
an integer program $\min \{ c^Tx \, : \, x \in \cF \}$. Suppose
that $x^* \notin \Z^d$ is an extreme point optimal solution to some 
convex relaxation to this
problem.  Then one attempts to cut-off $x^*$ by applying the following 
procedure.\\

\noindent First, a set $X = X(x^*) \subset \R^d$ with $x^* \in \interior(X)$ and $\interior(X) \cap \Z^d = \emptyset$ is identified.  Thus, denoting
$$ P \ \doteq \ X \cup \{ x \in \R^d \, : \, A x \not \ge  b \}$$
we have that $\cF \ \subseteq \ \R^d - \interior(P)$.  A valid inequality $\pi^T x \ge \pi_0$ is then sought such that
\begin{eqnarray}
(1) && \pi^T x = \pi_0 \ \ \mbox{supports} \ \conv(\R^d - \interior(P)), \ \ \mbox{and} \nonumber \\ 
(2) && \pi^T x^* < \pi_0.  \nonumber
\end{eqnarray}
One of the earliest versions of this
idea is embodied by the family of split cuts (see \cite{nemwols}), where the
set $X$ is the region bounded by two parallel hyperplanes. 
Intersection cuts \cite{bal71} are a closely related methodology.  For further material see \cite{andetal}, \cite{cornmarg}.
The disjunctive method for mixed-integer programming
(see \cite{bal75}, \cite{bal79}, \cite{balcercor93}), and the 
mixed-integer rounding procedure \cite{nemwols} rely on a similar
paradigm (also see
\cite{lacilex}, \cite{sheraliadams}).
We 
note that in the standard form of the 
above procedure integrality of the variables is only used in
the construction of the lattice-free set $X$; however the cut  
$\pi^T x \ge \pi_0$ is
computed simply using the geometry of the set $\R^d - \interior(P)$ (or
even, just the set $\R^d - X$) and
the point $x^*$.  Another point is that the set $X$ is typically quite simple
(e.g. a polyhedron defined by a small number of inequalities). Our approach to nonlinear, 
mixed-integer programs seeks to adapt the lattice-free approach, leading
to the study of sets $S$ as outlined above from a cutting-plane perspective.
In Section \ref{cardconst} we will describe a specific rendition of our
method in the context of cardinality-constrained convex
quadratic programming.\\

We apply the procedure closely related to \textit{lifting} valid inequalities for 
mixed-integer programs.  See
\cite{nemwols} for background.  It can be summarized
as follows. Let $\cF \subseteq \R_+^n$ be the feasible
region  of an integer program, and for $k < n$ let
$\sum_{j = 1}^k\beta_j x_j \le \beta_0$ be valid for $\cF$ when $ x_j = 0$ for all $j > k$.  Lifting is the process whereby this inequality is modified so as to yield an inequality $\sum_{j = 1}^{k+1}\beta_j x_j \le \beta_0$ valid for $\cF$ when $x_j = 0$ for all $j > k+1$.
Geometrically, 
the hyperplane defined
by $\sum_{j}^k \beta_j x = \beta_0$ is being rotated so as to support
$\conv(\cF)$ at a point $\hat x$ with $\hat x_{k+1} > 0$.  
 From this perspective, lifting may more
aptly be referred to as ``tilting'' (see the discussion in \cite{lifttilt}). 
Lifting techniques have proved compelling in that they are supported by
strong theory and can also provide a computationally
practicable way to strengthen valid inequalities.  As a result lifting is ubiquitous in mixed-integer programming solvers.

The extension of lifting to
the nonlinear or continuous
setting is not new: see \cite{danoq}, \cite{ismael1}, \cite{ismael2}, and \cite{alper}. Also see \cite{belmillnam}, which lifts ``tangent'' inequalities to approximate multilinear functions.   An interesting 
use of lifting appears in \cite{qbmarg}, which approximates, using
lifted linear inequalities, SDP relaxations of quadratically constrained sets.  
A comprehensive framework for coordinate-wise lifting 
is presented in \cite{richtaw}, where  given an arbitrary function $f(x,y) : \R^{p + n} \rightarrow \R$, and a linear underestimator for $f$ which is valid when
$y = 0$, i.e. $ f(x, 0) \, \ge \, \bar \alpha^T x  - \delta,$
lifting is used to modify this inequality 
so as to make
it valid for all $(x,y)$.  This results in an inequality of the form
$f(x,y) \, \ge \, \bar \alpha^T x  + \nu^Ty - \delta$, for an appropriate
vector $\nu$.  The authors discuss sequence independent lifting (over the
variables $y$) and applications to several problem types, such as bilinear
knapsack sets. \\

The techniques in this paper involve linear approximation to nonlinear functions
and further we focus on quadratics.  Both subjects have received significant
attention in the literature. One of the earliest results (see \cite{mccormick} and \cite{alkha})
is the characterization of the convex envelope of a box-constrained bilinear form $x_1x_2$ (also see \cite{rikun}).   These results lead to
techniques that have been incorporated in software systems such as BARON \cite{baron} and Couenne \cite{couenne}. 

The problems we consider fall within the broader scope of global optimization
problems. The work in \cite{talsah1} and \cite{talsah2} has resulted in 
key advances that can be applied in very general settings. An important
idea is that of outer approximation.  Given a convex 
function $g \, : \R \rightarrow \R$, \cite{talsah1} shows how to
construct an outer approximation to $g$ 
with arbitrary accuracy $\epsilon$ and using a number of lines that grows
(asymptotically) as $\epsilon^{-1}$.  Other techniques include the development
of effective relaxations to typical nonlinear functions (such as exponential,
logarithmic, and multilinear functions), the automatic generation of 
convex underestimators of general functions, and the adaptation of traditional
branch-and-bound to continuous domains.  These techniques are
amenable to implementation in a branch-and-cut setting (see \cite{talsah2}); 
 and have been included in BARON, achieving computational success.  Convex
extensions of a given function $f$ are considered in \cite{talsah3}; these are convex functions that agree
with $f$ on a subset of its domain.  This theory is
further used in \cite{talsah3} to study the convex envelope of the function $x/y$ of
two real variables $x, y$ over a rectangle in $\R^2$; additionally several results
are presented concerning convex envelopes of multilinear functions.

Recently, some interesting new 
results on multilinear forms have been obtained, see for example \cite{lnl}.  A survey 
is provided in \cite{burlet}; with additional material of interest 
in \cite{lee}, \cite{sax1}, \cite{sax2} and \cite{mitchelletal}.  A polyhedral approximation scheme
for nonconvex quadratically constrained quadratic programs is given in \cite{bao}.  Also see \cite{kojtunc1}.
A different,
frequently-applied construct is the Reformulation-Linearization Technique (RLT)
and semidefinite programming extensions;  see \cite{SA1} and \cite{SA2}.
Earlier work \cite{stubbsmehrotra}, \cite{frangent}, \cite{gunderoth} 
(also see \cite{ceriasoares}) in nonlinear 0-1 mixed-integer programming produced 
families of cuts arising from first-order
information.
Finally, the
connection with semidefinite programming has yielded a number of deep 
results focusing on quadratic functions, see for example \cite{kurtalone}, \cite{burerletchford}, \cite{anstreicherburer}. This provides
an alternative (but related) methodology for addressing 
some of the problems we consider.  

\subsubsection{Organization}
This paper is organized as follows.  In Section \ref{convexclosed} we study
the set $S$ defined in (\ref{Sdef}); we formally define our lifted 
inequalities and prove Theorem I given above. Then
we strengthen our result when $Q(x)$ satisfies an appropriate generalization
of strong convexity (Section \ref{psistrong}, obtaining Theorems II and III) with additional strengthening
when $Q(x)$ is a positive-definite quadratic in Section \ref{quadspecial}.
In Sections \ref{polyhedron} and \ref{ellipsoidal} we obtain Theorem IV (i);
presenting a polynomial-time separation algorithm
for $\conv(S)$ when $Q(x)$ is positive-definite quadratic and the set $P$
in (\ref{Sdef}) is a polyhedron and an ellipsoid, respectively.  Section 
\ref{twoquad} obtains a polynomial-time separation algorithm 
for a set $\{(x, w, q) \in \R^d \times \R \times \R \, : \, q \ge x^THx + h^Tx, \, w \le x^TAx \}$, where $H \succ 0$ and $A \succeq 0$, which is also a 
special case of (\ref{Sdef}) (Theorem IV (ii)).  Finally, in Section \ref{experiments} we present
initial experimental results.\\

\noindent We will use the following terminology: 

\begin{DE} \hspace{.1in}\\ 

\noindent (1) For a set $V \subseteq \R^d$, its {\em boundary} is defined as
$$\partial V = \{ x \in \R^d \, : \, 
\cB \cap V \neq \emptyset \ \mbox{and} \ \cB \cap (\R^d - V) \neq \emptyset \ \mbox{for every open ball $\cB$ containing $x$} \}.$$  

\noindent (2) Given a set $X \subseteq \R^d$, an inequality $\gamma^T x \ge \gamma_0$ which is valid for $X$ will be called a {\em supporting inequality for $X$} if $\{x \in \R^d \, : \, \gamma^Tx = \gamma_0 \}$ defines a supporting
hyperplane for $\conv(X)$. We will also say that $\gamma^T x \ge \gamma_0$ 
{\em supports} $X$, and if $y \in X$ is such that $\gamma^T y = \gamma_0$ 
we will also say that $\gamma^T x \ge \gamma_0$  supports $X$ {\em at} $y$. 
\end{DE}

\section{Lifted first-order cuts}\label{convexclosed}
Here we consider the set $S$ given by (\ref{Sdef}) where 
$Q$ is convex and differentiable and $P \subset \R^d$.  We will prove
a more detailed version of Theorem I, given as Theorem \ref{domination}, below. 
This theorem provides a characterization of supporting hyperplanes for
$\conv(S)$, in particular singling out the lifted inequalities (\ref{firstlifted0}).   
We will first introduce these lifted inequalities and prove a series of 
results (Lemma \ref{lemma1.2new}, and Propositions \ref{twice} - \ref{lflat}) leading to  Theorem \ref{domination}. Following this material, in Section \ref{psistrong} we provide
two results that hold when $Q(x)$ grows faster than linearly in 
every direction.
First, we obtain Theorem II (which characterizes polynomial-time separation from $\conv(S)$). 
Second, assuming (additionally) that the boundary of $P$ is appropriately structured, 
we prove that the lifted inequalities define hyperplanes guaranteed to support $S$ at (at least) two
different points (Theorem III in the Introduction).  Finally, in Section \ref{quadspecial} we discuss the case where $Q(x)$ is a positive-definite quadratic,
which allows for a geometric characterization of the lifted inequalities.\\

\noindent We first provide a brief motivation for our approach. We are interested in strengthening the linearized inequality (\ref{unlifted}) at
a point $y \in \partial P$ by modifying it in the form
\begin{eqnarray}
&& q  \ \ge  \ Q(y) + \nabla Q(y)^T(x - y)  +  2 p^T(x - y) \label{firstlifted}
\end{eqnarray}
for some $p \in \R^d$.  Note that this constitutes a strengthening only in the half-plane $\{ x \in \R^d \, : \, 
p^T(x - y) > 0\}$.  And in order for this strengthening to be valid for $S$ we must also have that 
\begin{eqnarray}
 \{\, x \in \R^d \quad : \quad  Q(x) \, <  \, Q(y) + \nabla Q(y)^T(x - y)  +  2 p^T(x - y) \, \} & \subseteq & \interior(P). \label{xyp}
\end{eqnarray}
Since $Q$ is differentiable it follows that for any $r \in \R^d$ with $p^Tr > 0$, any $x$ of
the form $x = y + \lambda r$ will be such that $(x, Q(x))$ violates (\ref{firstlifted}) provided $\lambda > 0$ is small enough (and how small may depend on $r$).  That is to say, for any $r$ with $p^Tr > 0$ there 
exist 
points $x = y + \lambda r$ with $\lambda > 0$
with $x$ contained in the set in the left-hand side of (\ref{xyp}).
We
now make these notions precise.

\begin{DE}\label{flatpoint} Let $y \in \partial P$ and $p \in \R^d$ with $\norm{p} = 1$. We say
that $P$ is locally flat at $y$ with normal $p$, if for every  $r \in \R^d$
with $p^T r > 0$ there exists $\epsilon(r) > 0$ such that
$$ y \, + \, \delta \, r \ \in \ \interior(P) \quad \forall \ \delta \ \mbox{with} \  0 < \delta \le \epsilon(r).$$
\end{DE}
Intuitively, $P$ being locally flat at $y$ with normal $p$ means that for any vector $r$ with 
positive inner product with $p$ we can move, starting at $y$, a positive distance ``into'' $P$ along $r$.  
\begin{EX}\label{locallyflatexample}
(a) If every connected component of $\partial P$ is a differentiable manifold
homeomorphic to $\R^{d-1}$ then $\partial P$ is locally flat at any
$y \in \partial P$, with a unique normal vector: the normal to the tangent space
to $\partial P$ at $y$, oriented into  $P$. (b) Let $P$ be a convex polygon in $\R^2$. Then $P$ is locally flat at every point on its boundary \emph{except} the vertices; using as normals the unit vectors 
normal to the facets, oriented into $P$. (c) The non-convex set $P = \{(x_1, x_2) \in \R^2 \ : \ \abs{x_1} \geq x_2\}$ is locally flat at every point on its boundary, even the vertex at $(0, 0)$ (with normal $(0,-1)$).  
\end{EX}

\noindent We can now begin our lifting construction.
\begin{DE} \label{hatalphadef0} Let $y \in \partial P$ be locally flat
with normal $p$.  For $\alpha \ge 0$ consider the 
inequality
\begin{eqnarray}
q & \ge & Q(y) + \nabla Q(y)^T(x - y) \, + \, 2 \alpha p^T(x - y) \label{partlifted}.
\end{eqnarray}
The lifting coefficient
at $y$, with respect to $p$, is $\hat \alpha =  \hat \alpha(P, p, y) \doteq  \sup \{ \, \alpha \, : \, \mbox{(\ref{partlifted}) \, \mbox{is valid for $S$}} \, \}.$ \end{DE}

\noindent We remark that one can equivalently 
write $\hat \alpha = \sup \{  \alpha \, : \, \mbox{(\ref{partlifted}) \mbox{is valid for $S$ for $0 \le \alpha \le \bar \alpha$}} \}$. Clearly the lifting coefficient is nonnegative, and we are interested in the
cases where this quantity is actually positive. 

\begin{EX}\label{ex1} Let $d = 2$, $Q(x) = x_1^2 + x_2^4$, and $P = \{ x \in \R^2 \, : \, x_2 \le |x_1| + 1 \ \mbox{and} \ x_2 \ge (x_1 - 1)(x_1 - 2)^2 \, \}$.
Thus $\nabla Q(y)^T = (2 y_1, 4 y_2^3)$  for any $y \in \R^2$.
For $v \in \R$, write $c(v) = (v - 1)(v - 2)^2$. Then $P$ is locally flat 
at $y = (1,0)^T$ with normal
$$\left( \frac{-1}{\sqrt{1 + [1/c'(1)]^2}} \, , \, \frac{1/c'(1)}{\sqrt{1 + [1/c'(1)]^2}}\right)^T \ = \ \left( \frac{-1}{\sqrt{2}} , \frac{1}{\sqrt{2}} \right)^T.$$
Hence (\ref{unlifted}) takes the form
$$ q \ \ge \ 1 \, + \, 2(x_1 - 1) \, + \, \hat \alpha \left( \frac{-1}{\sqrt{2}}(x_1 - 1) + \frac{1}{\sqrt{2}} x_2 \right).$$
A calculation shows that $\hat \alpha = \sqrt{2}$, and so the inequality 
can be rewritten as
$ q \ \ge \ x_1 \, + \, x_2,$ which is binding at $(x, Q(x))$ with $x = (1,0)^T$ and 
$x = (0,1)^T$.
\end{EX}

\begin{DE}\label{LFOineq} Suppose $P$ is locally flat at $y$ with normal $p$,
and that $\hat \alpha < + \infty$.  We call
\begin{eqnarray}
q & \ge & Q(y) \, + \, \nabla Q(y)^T(x - y) \, + \, 2 \hat \alpha p^T(x - y) \label{lifted}
\end{eqnarray}
a {\em lifted first-order (LFO) inequality} generated at $y$ (with respect to
$p$). 
\end{DE}

\noindent The following result establishes simple properties of
the lifting coefficient.

\begin{LE}\label{lemma1.2new} Let $P$ be locally flat at $y$ with normal $p$. {\bf(a)} If there exists $v \notin P$ such that $p^T(v - y) > 0$ then $\hat \alpha < +\infty$.
{\bf(b)} If $p^T(v - y) \le 0$ for all $v \notin P$ then $\hat \alpha = + \infty$.
{\bf(c)} If $\hat \alpha < +\infty$ then (\ref{lifted}) is valid for $S$.
\end{LE}
\noindent {\em Proof.}  {\bf(a)} Let $v \notin P$ satisfy $p^T(v - y) > 0$. 
Then clearly for $\alpha$ large enough, $(v, Q(v))$ will violate (\ref{partlifted}), and so $\hat \alpha < +\infty$. 
\noindent {\bf (b)}  This follows trivially since for any $\alpha > 0$ and
any $x \notin P$ the right-hand side of (\ref{partlifted}) is 
dominated by that of (\ref{unlifted}). {\bf(c)} This fact follows by continuity. \hspace{.1in}  \QED\\

\noindent {\bf Remarks.} The construction that culminates in Definition \ref{LFOineq} is a generalization of the classical lifting construction in 
mixed-integer programming (see \cite{nemwols}, \cite{ismael1}, \cite{ismael2}, 
\cite{alper}).  The LFO inequality at $y$ uses the local structure of
$P$ 
to strengthen the linearization inequality (\ref{unlifted}); the strengthening
is only local, however the LFO inequality is (globally) valid.  \\

\noindent We next prove the main result in this section. 
\begin{THM} \label{domination} Let
$\delta q  \ge  \beta^T x + \beta_0$ be valid for $S$ and 
binding at $(y, z) \in S$ for some $y \in \R^d - \interior(P)$ and $z \ge Q(y)$.  Then  at least one of the following conditions
holds: 
\begin{itemize}
\item [(1)] $\delta = 0$ and $ \beta^T x + \beta_0 \le 0$ is 
valid for $\R^d - P$. 
\item[(2)] $\delta > 0$, $z = Q(y)$ and $\delta q  \ge  \beta^T x + \beta_0$ is a positive multiple of the linearization inequality at $y$. 
\item[(3)] $\delta > 0$, $z = Q(y)$ and $y \in \partial P$. Moreover, $\delta q  \ge  \beta^T x + \beta_0$ is a nonnegative linear combination of the linearization inequality at $y$
and a linear inequality supporting $\R^d - \interior(P)$ at $y$. 
\item[(4)] $\delta > 0$, $z = Q(y)$, $y \in \partial P$ and $P$ is locally flat at $y$ with some normal $p$. Moreover
$\delta q \ge \beta^T x + \beta_0$ is implied by the linearization inequality at $y$ together with
the LFO inequality at $y$ with respect to $p$. 
\end{itemize}
\end{THM}

\noindent The proof of this theorem will be broken into a sequence of steps.
We will consider a fixed inequality
\begin{eqnarray}
\delta q & \ge & \beta^T x \, + \, \beta_0 \label{genvalid}
\end{eqnarray} 
that is supporting for $S$ at a point $(y, z)$ 
with $y \in \R^d - \interior(P)$ and $z \ge Q(y)$
and obtain Theorem \ref{domination} through
a sequence of results.  We begin with  some simple preliminary observations.\\

\noindent {\bf 1.} Since (\ref{genvalid}) is valid for $S$, we have $\delta \ge 0$.  If $\delta = 0$ then $\beta^T x  +  \beta_0 \le 0$ is valid for $\R^d - P$ and we are done (case (1) of Theorem \ref{domination}).  Thus we will assume $\delta > 0$ and by scaling if necessary that
$\delta = 1$. Hence we must have $z = Q(y)$ and in summary
 (\ref{genvalid}) is binding at $(y, Q(y))$.  \\

\noindent {\bf 2.} Write 
\begin{eqnarray}
 \beta^T x + \beta_0 & = & Q(y) + \nabla Q(y)^T(x - y) \, + \, 2 v^Tx \, - \, 2 v_0, \nonumber
\end{eqnarray}
for appropriate $v\in \R^d$ and $v_0 \in \R$.  Since (\ref{genvalid}) holds with
equality at $(y, Q(y))$ it follows that $ v_0 \, = \, v^T y$ and 
we can rewrite (\ref{genvalid}) as
\begin{eqnarray}
q & \ge & Q(y) + \nabla Q(y)^T(x - y) + 2 v^T(x - y). \label{adjusted}    
\end{eqnarray}

\noindent {\bf 3.} If $v = 0$ in (\ref{adjusted}) then (\ref{adjusted}) is the linearization inequality at $y$ (case (2) of Theorem \ref{domination}).   We will therefore assume $v \neq 0$. \\

\noindent We will use form (\ref{adjusted}) of (\ref{genvalid}), with $v \neq 0$, in Propositions \ref{twice}, \ref{partial}, \ref{lflat} and Corollary \ref{grand} given next.

\begin{PR} \label{twice} Let $w \in \R^d$ be arbitrary such that $v^T w > 0$.  Then 
there exists a positive value $\epsilon = \epsilon(w)$ such that for any $0 < \delta < \epsilon$ the right-hand side of (\ref{adjusted}) evaluated at 
$x = y + \delta w$ exceeds $Q(y + \delta w)$ since $v^Tv > 0$. \end{PR}
\noindent {\em Proof.} For $\delta > 0$ write $F(\delta) = Q(y + \delta w)$.  The right-hand side of (\ref{adjusted}) evaluated at 
$y + \delta w$ equals
$$Q(y) + \delta (\nabla Q(y) + 2v)^Tw \ \doteq \ G(\delta).$$
Then $G(0) = Q(y) = F(0)$ and $G'(0) = (\nabla Q(y) + 2v)^Tw > 
\nabla Q(y)^T w = F'(0)$.  As a result $F(\delta) < G(\delta)$ 
for $\delta > 0$ small enough. \QED

\begin{PR} \label{partial} We have $y \in \partial P$. \end{PR}
\noindent {\em Proof.}  Suppose by contradiction that $y \in \interior(\R^d - P)$. By Proposition \ref{twice}, for $\delta > 0$ small enough
the right-hand side of (\ref{adjusted}) evaluated at 
$y + \delta v$ exceeds $Q(y + \delta v)$.  This is a contradiction since
for 
$\delta > 0$ small enough
$y + \delta v \in \interior(\R^d - P) \subseteq \R^d - \interior(P)$. \QED

\begin{PR} \label{lflat} $P$ is locally flat at $y$, with normal $p \, \doteq \, \frac{v}{\norm{v}}$. \end{PR}
\noindent {\em Proof.} Follows directly from Propositions \ref{twice} and \ref{partial}.  \QED  \\

\noindent As a result of Proposition \ref{lflat}, we can rewrite (\ref{adjusted}) as 
\begin{eqnarray}
q & \ge & Q(y) + \nabla Q(y)^T(x - y) + 2 \norm{v} p^T(x - y). \label{adjusted2}    
\end{eqnarray}
where $p$ is as in Proposition \ref{lflat}.  Write 
$\hat \alpha \doteq \hat \alpha(P,p,y)$. 

\begin{CO}\label{grand} $\hat \alpha \ge \norm{v} > 0$. Furthermore
\begin{itemize}
\item[(a)] If $\hat \alpha = + \infty$ then $p^T(x - y) \le 0$ is valid for $\R^d - \interior(P)$ and (\ref{adjusted2}) is a nonnegative linear combination
of the linearization inequality at $y$, and $p^T(x - y) \le 0$.
\item[(b)] 
If $\hat \alpha = \norm{v}$ the LFO inequality at $y$ with respect to $p$
\begin{eqnarray}
q & \ge & Q(y) + \nabla Q(y)^T(x - y) + 2 \hat \alpha p^T(x - y) \label{liftedhohoho}
\end{eqnarray}
and (\ref{adjusted2}) are identical constraints. 
\item[(c)] Suppose
$\norm{v} < \hat \alpha < +\infty $. Then at any $x$ with $p^T(x - y)> 0$ the
right-hand side of (\ref{liftedhohoho})
is strictly larger than that of (\ref{adjusted2}). Further, at
any $x$ with $p^T(x - y) < 0$ the right-hand side of the linearization inequality
at $y$ is strictly larger than that of (\ref{adjusted2}).
\item [(d)] Suppose
$\norm{v} < \hat \alpha < +\infty$. For any $x \in \R^d$
inequality (\ref{adjusted2}) is weaker than the combination of the linearization
inequality at $y$ and inequality
(\ref{liftedhohoho}). 
\end{itemize} 
\end{CO}
\noindent {\em Proof.} The definition of lifting
coefficient and validity of (\ref{adjusted2}) implies $\hat \alpha \ge \norm{v}$.  (a) Part (a) of Lemma \ref{lemma1.2new} implies that 
$0 \ge p^T(x - y)$ is valid for $\R^d - \interior(P)$; the rest of the
statement follows since clearly (\ref{adjusted2}) is obtained by adding $\norm{v}$ times $0 \ge p^T(x - y)$ to the linearization inequality at $y$. Part (b) is clear, and (c) follows since $\hat \alpha > 0$. Finally, (d) is a corollary of
(c): given $x \in \R^d$, if $p^T(x-y) \le 0$ the right-hand side of (\ref{adjusted2}) evaluated
at $x$, is at most that of the the linearization inequality at $y$.
And if $p^T(x -y) > 0$ then the right-hand side of (\ref{adjusted2}) evaluated
at $x$ is less than that of (\ref{liftedhohoho}). \QED \\

\noindent Note that (a) of Corollary \ref{grand} amounts to case (3) of Theorem \ref{domination},
and (b) and (c) to case (4). Thus we have completed the proof of Theorem \ref{domination}. 
\begin{RE} \label{super0} There are a number of conditions under which the
separation problem for $\conv(S)$ is equivalent to separation by LFO inequalities,
linearization inequalities, and valid inequalities for $\R^d - P$.  We will return to
this issue in Theorem \ref{superbasic}, below.\end{RE}

\begin{EX}\label{ex2} As an illustration of the use of LFO inequalities, consider Example \ref{ex1}.  Suppose we apply the
following heuristic for the problem $\min \{ Q(x) \, : \, x \in \R^2 - \interior(P) \}$. We start with the relaxation 
$$ \min \{ q \, : \, q \ge x_1 + x_2, \ q \ge 0, \ x \in \R^2\}$$
consisting of the LFO inequality at $(1,0)^T$ and the linearization 
inequality at $(0,0)$.  We initialize $\breve x_1 = 1$.  Then
we perform the following steps. \\
\noindent {\bf 1.} Solve the relaxation, obtaining solution $x^*$.\\
\noindent {\bf 2.} {\bf If} $x^* \in \R^2 - \interior(P)$, add to the 
relaxation the linearization inequality at $x^*$, and go to 1. {\bf Otherwise:}\\
\noindent {\bf 3.} Update $\breve x_1 = \frac{1}{2}(x^*_1 + \breve x_1)$.\\
\noindent {\bf 4.} Add to the relaxation the LFO inequality and the linearization  inequality at $(\breve x_1, (\breve x_1 - 1)(\breve x_1 - 2)^2)^T$.\\
\noindent {\bf 5.} Go to 1. \\
This heuristic will produce the following sequence of values $\breve x$ (truncated to three digits):\\ 
\noindent $1.000, 0.500, 0.545, 0.596, 0.648, 0.702, 0.687, 0.691, 0.694, 0.696$.  After nine iterations, the formulation proves a lower bound of
$0.55577$ on the value of the optimization problem.  Moreover, setting
$\hat x_1 = 0.696208$ and $\hat x_2 = (\hat x_1 - 1)(\hat x_1 - 2)^2$, we
have that $\hat x$ is feasible while $Q(\hat x) \approx  0.55582$.
\end{EX}

\subsubsection{Strong convexity implications on LFO inequalities}\label{psistrong}
Here we address two issues that arise from the above analysis
and which are resolved when $Q(x)$ satisfies
a strong convexity assumption: the relationship between separation from LFO inequalities and separation from $\conv(S)$ (Theorem II in the Introduction), and whether LFO inequalities are binding at more 
than one point (Theorem III).  To formalize our approach we first 
review some standard concepts. \\

\noindent A function $f \, : \, \R^d \rightarrow \R$ is called 
strongly convex (with modulus $2$) \cite{boydvan} if
$$f(x) \ \ge \ f(y) + \nabla f(y)^T(x - y) + \norm{x - y}^2, \quad \forall x, y \in \R^d.$$
Strongly convex functions are of interest because they include positive-definite
quadratics, which we will focus on in some of our results, below.  A generalization
of strong convexity is $\psi$-strong convexity, where for $\psi \, : \, \R_+ \rightarrow \R$, the following condition is satisfied:
$$f(x) \ \ge \ f(y) + \nabla f(y)^T(x - y) + \psi(\norm{x - y}), \quad \forall x, y \in \R^d.$$
This is a generalization because a strongly convex function is 
$\psi$-strongly convex with $\psi(t) = t^2$, but
for example, for $d = 1$, $f(x) = x^{4}$ is not strongly convex (near $x = 0$) 
but is $\psi$-strongly convex with $\psi(t) = t^4$. See \cite{alphastrong} 
for a discussion of generalized strong convexity.    Here, we will be relying on $\psi$-strongly convexity where 
$\psi$ is a function satisfying the following conditions:
\begin{eqnarray}
&& \psi \ \mbox{is strictly increasing,} \quad \quad \psi(0) = 0, \quad \quad \mbox{and} \  \lim_{t \rightarrow +\infty} \frac{\psi(t)}{t} = +\infty, \label{growsfast}
\end{eqnarray}
Under this criterion, $\psi$-strongly convex functions grow faster than linearly in every direction. \\

\noindent We now turn to the first issue raised above.  
Whereas Theorem \ref{domination} classifies supporting
inequalities for the set $S$ it does not directly address separation from
$\conv(S)$.   The next result and corollary address this issue.

\begin{THM}\label{superbasic} 
Suppose $Q(x)$ is $\psi$-strongly convex where $\psi$ satisfies (\ref{growsfast}), and
that $(\bar x, \bar q) \in \R^d \times \R$ satisfies the following conditions: 
\begin{itemize}
\item [(a)] $(\bar x, \bar q) \notin \closure(\conv(S))$.
\item [(b)] $\bar x \in \conv( \R^d - \interior(P))$.
\item [(c)] $(\bar x, \bar q)$ satisfies the linearization inequality at 
$\bar x$.
\end{itemize}
Then there is an LFO inequality that separates 
$(\bar x, \bar q)$ from $\conv(S)$.\end{THM}
\noindent {\em Proof.} By (a) $(\bar x, \bar q)$
violates an inequality $\delta q \ge \beta^T x + \beta_0$
which is valid for $conv(S)$.   By (b), $\delta > 0$ and without
loss of generality $\delta = 1$.  Choose $z \in \R^d - \interior(P)$; thus
$Q(z) \ge \beta^T x + \beta_0 + \epsilon_0$ for some $\epsilon_0 \ge 0$.  Furthermore,
since $\psi$ satisfies (\ref{growsfast}), we have that there exists $R \ge \norm{z}$ such that
$$ Q(x) \ \ge \ \beta^T x + \beta_0 + \epsilon_0, \quad \quad \forall \ x \ \mbox{with} \ \norm{x} \ge R.$$
Consequently, there exists $y \in \R^d - \interior(P)$ with $\norm{y} \le R$ and $0 \le \epsilon_1 \le \epsilon_0$ such that
\begin{eqnarray}
Q(y) & = & \beta^T y + \beta_0 + \epsilon_1 \quad \mbox{and} \label{itworked} \\
 q & \ge & \beta^T x + \beta_0 + \epsilon_1  \quad \mbox{is valid for $S$.} \nonumber
\end{eqnarray}
We will now classify (\ref{itworked}) as per Theorem \ref{domination}.
Clearly case (1) of Theorem \ref{domination} does not apply. Further,
by convexity of $Q(x)$, if the linearization inequality at any $z \in \R^d - \interior(P)$  is violated by $(\bar x, \bar q)$, then so is the linearization inequality at $\bar x$ itself, a contradiction by assumption 
(c).  Thus cases (2) and (3) of Theorem \ref{domination} do not apply.  We conclude as desired. \QED 

\begin{CO} \label{super} Suppose $Q(x)$ is $\psi$-strongly convex where $\psi$ satisfies (\ref{growsfast}), 
that there is a polynomial-time
separation oracle for $\R^d - \interior(P)$, and that $\nabla Q(x)$ is polynomial-time computable at any
$x$.  Then polynomial-time separation over $\conv(S)$ is
equivalent to polynomial-time separation over the LFO inequalities. \end{CO}
\noindent {\em Proof.} Given point $(\bar x , \bar q) \in \R^d \times \R$ we can check 
in polynomial time whether it satisfies the
linearization inequality at $\bar x$ as well as $\bar x \in \conv(\R^d - P)$.  \QED\\

\noindent We now address the second issue raised above,
where a $\psi$-strong convexity assumption again has
a significant implication.
A question left open in Section \ref{convexclosed} is whether,
given
an LFO inequality obtained at a point $y$, there exists
$w \neq y$ and feasible such that the lifted inequality is also binding at
$(w, Q(w))$.
In Section \ref{polyhedron} we will see a specific case where the fact
that this condition holds is used to construct a polynomial-time 
separation algorithm for $\conv(S)$.  However, and in contrast to what
happens in the linear mixed-integer programming setting, the following two examples show that 
the condition does not always hold, for two possible reasons illustrated in 
the following examples.

\begin{EX}\label{notattained}  Let $P = \{ x \in \R^2 \, : \, -x_1^2  \le x_2 \le 1 + e^{-x_1} \},$ and $Q(x) = x_2 + e^{-x_2} -1$.  Then $0 \in \partial P$, $P$ is locally flat at $0$ with normal $(0,1)^T$, $Q(0) = 0$ and $\nabla Q(0) = 0$.
Thus the LFO inequality at $0$ has the form
$q \ge \hat \alpha x_2$.  
Furthermore when $x_2 >1$, $Q(x) > e^{-1}x_2$.  It follows that $\hat \alpha = e^{-1}$ since with 
this choice $q \ge \hat \alpha x_2$ is valid, but any larger value will exclude points in $S$.  However
any point $x \neq 0$ with $Q(x) = e^{-1}x_2$ satisfies $x_2 = 1$ and thus $x \in \interior(P)$.
\end{EX}

\begin{EX}\label{notflatenough} Let $P = \{x \in \R^2 \, : \, x_1 \ge 1 \} \cup \{x \in \R^2 \, : \, 0 \le x_1 \le 1 \ \mbox{and} \ |x_2| \le (2x_1 - x_1^2)^{1/2} + x_1 \}$.  
Suppose $Q(x) = \norm{x}^2$. By construction, the ball with center $(1,0)^T$ 
and unit radius is contained in $P$, and so $P$ is locally flat at $0$, with
unique normal $(1,0)^T$, and we obtain the LFO inequality $q \, \ge \, \hat \alpha x_1$. \\

\noindent Now, for any $R > 0$, the inequality $q \ge 2 Rx_1$ is valid for $S$ if and
only if $x \in \interior(P)$ whenever $\norm{x}^2 < 2 R x_1$,  i.e. whenever the ball with center $(R, 0)^T$ and radius $R$ is
contained in $P$. Therefore $\hat \alpha \ge 1$. However, for any value $R > 1$
we can find $0 < x_1 < 1$ such that
$$ (2R x_1 - x_1^2)^{1/2} > (2x_1 - x_1^2)^{1/2} + x_1.$$
It follows that for \textit{any} $R > 1$ the ball with radius $R$ and center
at $(R,0)^T$ is not contained in $P$.  As a result $\hat \alpha = 1$ and
yet the only point $(x, Q(x))$ with $x \in \R^d - \interior(P)$ where 
$q \ge 2 x_1$ 
 is binding is $(0,0)$.
\end{EX}

In Example \ref{notattained}, the function $Q(x)$ 
effectively grows at a linear rate in $x_2$,  while in Example \ref{notflatenough} the boundary
of $P$, near $0$, is curved too steeply in the direction of the lifting  
(so that for $R = 1 + \epsilon$ with $\epsilon > 0$ and small, the inequality
$q > 2 R x_1$ is violated by points whose norm tends to zero as $\epsilon \rightarrow 0$).   

We will show
next that if (1) $Q(x)$ grows faster than linearly, in every direction, and (2) $\interior(P)$ contains a half-ball with
positive radius and center
at each point where lifting is performed, then any LFO inequality is binding at (at least) two different points,
obtaining Theorem IV of the Introduction.  This will be done in Theorem \ref{thm1.2strong} below. We first prove a simple technical result concerning 
the lifting construction.  Given a function $\psi$ satisfying (\ref{growsfast}), for $k > 0$ we define 
$$\chi^{\psi}(k) \ \doteq \ \sup\{ t \ge 0 \, : \, \psi(t) < 2 k t \}$$
which is finite by (\ref{growsfast}).
\begin{LE} \label{lemma0new} Assume $Q(x)$ is $\psi$-strongly convex with 
$\psi$ satisfying (\ref{growsfast}).
Let $y \in \R^d$ and $0 \neq p \in \R^d$. Suppose
$x \in \R^d$ satisfies
\begin{eqnarray}
Q(x) & < & Q(y) + \nabla Q(y)^T(x - y) \, + \, 2 p^T(x - y). \nonumber
\end{eqnarray}
Then $\norm{x - y} \le \chi^{\psi}(\norm{p})$. \end{LE}
\noindent {\em Proof.} 
Using $\psi$-strong convexity
yields 
\begin{eqnarray}
Q(y) + \nabla Q(y)^T (x - y) +  \psi(\norm{x - y}) & < & Q(y) + \nabla Q(y)^T(x - y) \, + \, 2 p^T(x - y),  \quad \mbox{i.e.} \nonumber \\
\psi(\norm{x - y}) & < & 2 p^T(x - y) \ \le \ 2 \norm{p} \norm{x - y}. \nonumber
\end{eqnarray}
\QED\\

\noindent  Below we will use the following notation: given vectors $v \, \mbox{and nonzero} \ p \in \R^d$, and a real $\gamma > 0$, write 
$$ \cH(v, \gamma) \, \doteq \, \{ x \in \R^d \, : \, p^T(x - v) > 0, \ \norm{x - v} < \gamma \},$$
which is an (open) half-ball with center $v$ and radius $\gamma$.
\begin{THM}\label{thm1.2strong} Suppose $Q(x)$ is $\psi$-strongly convex where $\psi$ satisfies (\ref{growsfast}).  Let
$y \in \partial{P}$ and $p \in \R^d$
satisfy:
\begin{itemize}
\item[(i)] There exists $v \notin P$ such that $p^T(v - y) > 0$.
\item[(ii)] There exists a real $\gamma > 0$ such that $\cH(y, \gamma)  \subseteq P$.
\end{itemize}
Then $P$ is locally flat at $y$ with normal $p$, and
\begin{itemize}
\item[(1)] there exists $w \notin \interior(P)$ with $w \neq y$ and such that the LFO inequality at $y$
is binding at $(w, Q(w))$, i.e.
$$ Q(w) \ = \ Q(y) + \nabla Q(y)^T(w - y) \, + \, 2 \hat \alpha p^T(w - y).$$
\item[(2)] Furthermore, $\hat \alpha > 0$ and $(w - y)^Tp > 0.$
\end{itemize}
\end{THM}
\noindent {\em Pproof.}  Condition (ii) implies that $P$ is locally flat at
$y$. Furthermore, 
Lemma \ref{lemma1.2new} (b) implies that $\hat \alpha < + \infty$. 
But by definition of $\hat \alpha$, for any $\epsilon > 0$ there exists
$x^{\epsilon} \notin \interior(P)$ such that
\begin{eqnarray}
Q(x^{\epsilon}) & < & Q(y) + \nabla Q(y)^T(x^{\epsilon} - y) \, + \, 2 (\hat \alpha + \epsilon) p^T(x^{\epsilon} - y). \label{epsilonlevel}
\end{eqnarray}
Using Lemma \ref{lemma0new} we obtain $\norm{x^{\epsilon} - y} \le \chi^{\psi}(\hat \alpha + \epsilon) $.  Note that for any pair of values $0 < \epsilon < \delta$ we have $ \chi^{\psi}(\hat \alpha + \epsilon) \le \chi^{\psi}(\hat \alpha + \delta)$.
We conclude that since $\R^d - \interior(P)$ is closed, as $\epsilon \rightarrow 0$ there is an accumulation point $w \in \R^d - \interior(P)$ of the
points $x^{\epsilon}$.  Since $Q(x)$ is continuous, from (\ref{epsilonlevel})
we obtain
\begin{eqnarray}
&& Q(w) \ \le \ Q(y) + \nabla Q(y)^T(w - y) \, + \, 2 \hat \alpha p^T(w - y). \label{hereitsflat}
\end{eqnarray}
But since (\ref{partlifted}) is valid at $\alpha = \hat \alpha$, and we have
$w \in \R^d - \interior(P)$, we conclude that (\ref{hereitsflat}) holds as
an equality.  Further, by assumption (ii), $\norm{w - y} \ge \gamma > 0$, which implies $w \neq y$, and since $\psi$ is strictly increasing we have
(2) as well. \hspace{.1in}  \QED

\begin{RE} Condition (b) is a strengthening on the locally flat requirement
in Definition \ref{flatpoint}, which required that for each vector $r$ with 
with $p^T r > 0$ there exists $\epsilon(r) > 0$ such that $y + \delta r \in \interior(P)$ for all $0 < \delta < \epsilon(r)$.  In condition (b) we have $\epsilon(r) = \gamma$ for all appropriate $r$.  It is possible to relax condition (b) somewhat, to better account for the relationship between the values $\psi(\epsilon(r))$ and $p^T r$, for all $r$.
\end{RE}

\subsection{Specialization when Q(x) is a positive-definite quadratic}\label{quadspecial}

When $Q(x)$ is a positive-definite quadratic the constructions above can
be simplified and strengthened.  By changing 
coordinates if necessary we assume $Q(x) = \norm{x}^2$,  and thus
the LFO inequality at
a point $y$ with respect to a unit vector $p$ has the form $q \ge \norm{y}^2 + 2 (y + \hat \alpha p)^T(x - y)$. In this section we will show that $\hat \alpha > 0$ if and only if there exists a ball $\cB$ such that
$y \in \partial \cB$ and $\cB \subseteq P$ (Corollaries \ref{hee} and \ref{haw}, below). We
will also obtain other structural results, in particular a geometrical
interpretation of LFO inequalities. In Section \ref{polytime} we will use these improvements to 
obtain polynomial-time separation procedures for $\conv(S)$ when $P$ is 
either a polyhedron or an ellipsoid and $Q(x)$ is positive-definite quadratic.\\

\noindent We will use the following notation: given $\mu \in \R^d$ and $R \ge 0$, we  
write $\cB(\mu, R) \ = \ \{ \, x \in \R^d \, : \, \| x - \mu \| \le R \, \}$.  
The following property will be used in the sequel.

\begin{RE} \label{prop0prime} Let $y \in \R^d$ and $v \in \R^d$. Then
$x \in \R^d$ satisfies
\begin{eqnarray}
\norm{x}^2 & \le & \norm{y}^2 + 2y^T(x - y) \, + \, 2 v^T(x - y). \label{violated}
\end{eqnarray}
if and only if $x \in \cB(y + v, \norm{v})$ with equality in 
(\ref{violated}) iff $\norm{x -(y+v)} = \norm{v}$.
\end{RE}
\noindent {\em Proof.}  We can restate (\ref{violated}) as
$ \norm{x}^2 \, \le \, 2(y + v)^T x - \norm{y}^2 - 2 y^T v,$ from which the result follows. \QED \\

\noindent From this observation we obtain two corollaries:

\begin{CO} \label{hee} Suppose $P$ is locally flat at $y \in \partial P$ with normal $p$,
and that $\hat \alpha = \hat \alpha(P,p,y) > 0$ is finite.  Then:  
\begin{itemize}
\item[(a)]
$\cB(y + \hat \alpha p, \hat \alpha ) \subseteq P$. 
\end{itemize}
Further, suppose that there is a real $\gamma > 0$  such that $\cH(y, \gamma) \subseteq P$. Then
there exists $z \neq y$, $z \notin \interior(P)$ satisfying conditions (b) and
(c) given next.
\begin{itemize}
\item[(b)]  $\norm{z - (y + \hat \alpha p)} = \hat \alpha$, and 
\begin{eqnarray}
\norm{z}^2 & = & \norm{y}^2 \ + \ 2(y + \hat \alpha p)^T(z - y). \label{tightq}
\end{eqnarray}
Further, $z \in \partial P$, and writing $q = z - (y + \hat \alpha p)$, $P$ is locally flat at $z$, with normal $v = q/\norm{q}$.  
\item[(c)] $\hat \alpha(P, v, z) = \hat \alpha$.
\end{itemize}
\end{CO}
\noindent {\em Proof.} (a) This follows from Remark
\ref{prop0prime}.  (b) Let $z \in \R^d - \interior(P)$ be a vector that satisfies the LFO
inequality at $y$ with respect to $p$ with equality,  with $z \neq y$, which is guaranteed to 
exist by Theorem \ref{thm1.2strong}. Then by construction $z$ satisfies
(\ref{tightq}), and by Remark \ref{prop0prime},  $\norm{z - (y + \hat \alpha p)} = \hat \alpha$. Hence by part (a), $z \in \partial \cB(y + \hat \alpha p, \hat \alpha )$ and so $P$ is locally flat at $z$ with normal $v$. (c) By (a) and (b), 
$\hat \alpha(P, v, z) \ge \hat \alpha$. But any larger value of 
$\hat \alpha(P, v, z)$ would cut off $(y, \norm{y}^2)$. \QED

\begin{CO} \label{haw} Suppose $y \in \partial P$ and let $v \in \R^d - \{0\}$ be such that 
$\cB(y + v, \norm{v}) \subseteq P$.  Then $P$ is locally flat at $y$ with
normal $p = v/{\norm{v}}$ and $\hat \alpha(P,p,y) > 0$. 
\end{CO}
\noindent {\em Proof.} Since $\cB(y + v, \norm{v}) \subseteq P$ it clearly 
follows that $P$ is locally flat at $y$ with normal $p$.  Moreover, by
Remark \ref{prop0prime}, $\norm{x}^2 \, \ge \, \norm{y}^2 + 2y^T(x - y) \, + \, 2 v^T(x - y)$ holds for any $x \in \R^d - \interior(P)$.  Hence
$\hat \alpha(P, p, y) \ge 0$. \QED \\

The significance of Corollaries \ref{hee} and \ref{haw} is the following:  in the case that $Q(x) = \norm{x}^2$ an LFO inequality obtained at a point
$y \in \partial P$ has positive lifting coefficient if and only if there is a ball
$\cB(\mu, \sqrt{\rho}) \subset P$, with $\rho > 0$ and such that 
$\norm{y - \mu}^2  =  \rho$; furthermore the existence of such a ball 
implies
that  
$P$ is locally flat at $y$ with normal $(\mu - y)/\norm{\mu - y}$.    This is a sharpening of
Theorem \ref{domination} in that it simplifies the separation problem
for LFO inequalities.    For example, if $P$ is a 
polyhedron then we only need to consider LFO inequalities defined at points 
$y$ which are in the relative interior of some facet.

\vspace{.1in}
To conclude this section we point out a geometric characterization
of a set of the form $\R^d - \interior(P)$ which can be used to derive
valid inequalities for the corresponding set $S$ when $Q(x) = \norm{x}^2$.  Let $P$ be given. 
Then for any  $x \in \R^d$,
\begin{eqnarray}
&& x \in \R^d - \interior(P) \ \ \mbox{if and only if} \ \| x - \mu \|^2 \ge \rho, \ \ \mbox{for each ball $\cB(\mu, \sqrt \rho)$ contained in $P$.} \label{geometric}
\end{eqnarray}
As a result 
\begin{eqnarray}
&& q \, \ge \, 2 \mu^T x + \rho  - \| \mu \|^2, \ \ \mbox{for each ball $\cB(\mu, \sqrt \rho)$ contained in $P$} \label{geometric-q}
\end{eqnarray} 
is a family of inequalities valid for $S$.  We will term these the {\em ball inequalities}.  In fact, LFO inequalities are ball inequalities:
 consider the LFO inequality at a point $y \in \partial P$ with 
normal $p$ and lifting coefficient $\hat \alpha$.  Writing $\mu = y + \hat \alpha p$ and $\rho = \norm{ \mu - y}^2$, we have by Remark \ref{prop0prime}
that the inequality is violated 
precisely by those points $(x, \norm{x}^2)$ with $x$ in the interior
$B(\mu, \sqrt{\rho})$.  Thus the LFO inequality can be written
\begin{eqnarray}
&& q \ \ge \  2 \mu^T x \, + \, \rho \, - \, \norm{\mu}^2, \label{LFOball}
\end{eqnarray}
which
is a ball inequality.
Hence, even though (\ref{geometric}) provides a geometric characterization
of membership in $\R^d - \interior(P)$, Theorem \ref{domination} implies
that only a subset of the ball inequalities (\ref{geometric-q}) are needed 
to characterize $\conv(S)$: the
LFO inequalities.

\section{Polynomial-time separable cases}\label{polytime}
In this section we consider two cases where the 
characterization in Theorem \ref{domination} of
a set $\conv(S)$ with $S$ as in (\ref{Sdef}) leads to polynomial-time 
separation algorithms:  $Q(x)$ is positive-definite
quadratic and $P$ either a polyhedron or an ellipsoid (handled in Sections
\ref{polyhedron} 
and \ref{ellipsoidal}, respectively).

Here we recall the implications of Theorem \ref{superbasic} and Corollary \ref{super} concerning the separation problem.  Given $(x^*, q^*) \in \R^d \times \R$ 
we can trivially check whether this point satisfies all linearization
inequalities simply by checking if it satisfies the linearization 
inequality at $x^*$.  Assuming that $x^*$ also satisfies all valid
inequalities for $\R^d - P$, the only nontrivial case of the separation
problem is that where $x^* \in \interior(P)$, and we need to verify that $(x^*, q^*)$ satisfies all LFO inequalities.

\subsection{Polynomial-time separation of LFOs when P is a polyhedron and Q(x) is positive-definite quadratic}\label{polyhedron}
As before we can assume $Q(x) = \norm{x}^2$. Without loss of generality $\interior(P) \neq \emptyset$ and so 
$P$ is full-dimensional, so that $P \, = \, \{ \, x \in \R^d \, : \, a_i^T x \ge b_i, \ 1 \le i \le m \}$, where each inequality is facet-defining and nonredundant.  We assume $d \ge 2$ and $m \ge 2$.

\noindent For $1 \le i \le m$ let $\bar P^i = \{ \, x \in \R^d \, : \, a_i^T x \le b_i \}$; thus
$\R^d - P = \bigcup_i \bar P^i$.  Further, for $1 \le i \le m$ write:
$$ \bar Q^i = \{ \, (x, q) \in \R^d \times \R \, : \, a_i^T x \le b_i, \ q \ge \|x\|^2 \}.$$
Thus, $(x^*, q^*) \in \conv(S)$ if and only if $(x^*, q^*)$ can be written as a convex
combination of points in the sets $\bar Q^i$.  This is the disjunctive approach pioneered in 
Ceria and Soares \cite{ceriasoares}. Also see \cite{stubbsmehrotra}, \cite{frangent}, \cite{gunderoth}.  The resulting separation problem is carried out
by solving a second-order cone program with $m(d + 2)$ variables, $m$ conic constraints 
and  $d + 1$ linear constraints, and 
then using second-order cone duality in order to obtain a linear inequality (details in \cite{alex2012}).  In the specific context of this section, the work in this paper can be seen as a generalization of that in \cite{stubbsmehrotra}, \cite{frangent}, \cite{gunderoth}, all of which use the geometry of $0-1$ disjunctions and first-order estimates in order to derive cuts.

Here we will present an algorithm that, given $x^* \in \interior(P)$, finds an LFO inequality 
which proves the strongest lower bound on $q$ at $x^*$, that is to say, an LFO inequality 
$q \ge \beta^T x + \beta_0$ whose right-hand side is maximized at $x^*$. This requires solving $m-1$ 
convex quadratic programs, each with $d + 1$ variables and $2m - 1$ linear constraints.  The quadratic programs are given below in formulation SEP(i), for $1 \le i < m$.  The potential advantage of this approach relative to the use of the disjunctive formulation is twofold: computations are
done in the original space of variables, and the separation problem is of a 
simpler nature.  A numerical comparison between the two methods will 
be provided in Section \ref{disjunctcomparo}.\\

Our main construction is given in Lemma \ref{slice} below. In order to motivate our approach we first present some introductory remarks.
First, any LFO inequality is generated at some point $y \in \partial P$ where $P$ is locally flat (with some normal $p$). This property holds iff $y$ is in the relative interior of one of the facets defining $P$, say the facet 
corresponding to inequality $a_i^Tx \ge b_i$, in which case $p = a_i/\norm{a_i}$. 
Moreover, by Corollary \ref{hee} there is a ball contained in $P$, with radius
equal to the lifting coefficient, which contains in its boundary both $y$
and another point $z \in \R^d - \interior(P)$ .  
Using Corollaries \ref{hee} and \ref{haw}, necessarily we must then have
that $z$ is in the relative interior of another facet of $P$, say 
the facet defined by $a_j^Tx \ge b_j$, for some $j \neq i$.  And, moreover,
there is a symmetric relationship between $y$ and $z$, in the sense that 
lifting from $z$ will produce the same ball (as per Corollary \ref{hee}) and 
the same lifting coefficient.

These observations suggest that we explore the interaction between pairs of
facets.  To that effect, writing for any pair of distinct indices $1 \le i \le m$, $1 \le j \le m$
\begin{eqnarray}
&& P^{i,j} \ \doteq \ \{ x \in \R^d \, : \, a_i^T x \ge b_i, \ a_j^T x \ge b_j \}, \ \ 
\end{eqnarray}
we then have:
\begin{PR} \label{trivial} Let $1 \le i \le m$ and let $y \in P$ be in the relative interior of the facet defined
by $a_i^T x \ge b_i$.  Then $\hat \alpha(P, a_i/\norm{a_i}, y) \ = \ \min_{j \neq i} \hat \alpha(P^{i,j}, a_i/\norm{a_i}, y)$.
\end{PR}
{\em Proof sketch.} Any $0 \le \alpha \le \min_{j \neq i} \hat \alpha(P^{i,j}, a_i/\norm{a_i}, y)$ is a valid lifting coefficient at $y$ (with respect to the set 
$P$); and any larger coefficient will exclude at least one feasible point. \QED

We now use these
observations to obtain a characterization of the
LFO inequalities that leads to polynomial-time separability.  In particular, Lemma \ref{slice} below will allow us
to compute the minimum in Proposition \ref{trivial} (for a given $1 \le i \le m$)
by solving a convex quadratic program with $d+1$ variables and $m$ linear
constraints.  For $1 \le i \le m$ let 
$H_i  \doteq  \{ \, x \in \R^d \, : \, a_i^Tx = b_i \}$.  For
$i \neq j$ let $H_{\{i,j\}} \doteq \{ \, x \in \R^d \, : \, a_i^Tx = b_i, \  a_j^Tx = b_j\}$.

\begin{LE}\label{slice}  Let $1 \le i \le m$, $1 \le j \le m$ be distinct
and let $y$ be in the relative interior of the facet of $P$ defined by
$a_i^T x \ge b_i$.
Then there exists polynomially computable 
$p_{ij} \in \R^d$ and $q_{ij} \in \R$ such that
\begin{eqnarray}
\hat \alpha_{ij} \ \doteq \ \hat \alpha(P^{i,j}, a_i/\norm{a_i}, y) & = & p^T_{ij} y + q_{ij}. \nonumber
\end{eqnarray} 
Further, for any $v \in H_{\{i,j\}}$, $p^T_{ij} v + q_{ij} = 0$.
\end{LE}
\noindent {\em Proof.} By
Corollary \ref{hee} $\hat \alpha_{ij}$ is equal to the largest radius of a ball that can be inscribed in $P^{i,j}$, with $y$ in its
boundary. Thus, if $H_{\{i,j\}} = \emptyset$, i.e. $H_i$ and $H_j$ are parallel,
$\hat \alpha_{ij}$ equals
half the distance between $H_i$ and $H_j$, and is therefore independent of
$y$, i.e. it is trivially an affine function of $y$.  In what follows we assume 
$H_{\{i,j\}}\neq \emptyset$. Set
\begin{eqnarray}
&& \mu \doteq y + \hat \alpha_{ij} \frac{a_i}{\norm{a_i}},  \quad \mbox{and} \quad \rho \doteq \hat \alpha_{ij}^2 = \norm{\mu - y}^2.\label{mudef}
\end{eqnarray}
Then by Corollary \ref{hee}
\begin{eqnarray}
&& \cB(\mu, \sqrt{\rho}) \subseteq P^{i, j}, \quad  \cB(\mu, \sqrt{\rho}) \cap H_i = y, \quad \mbox{and} \quad \cB(\mu, \sqrt{\rho}) \cap H_j = z, \ \mbox{for some $z \in H_j$}. \label{important1}
\end{eqnarray} 
We have that $H_{\{i,j\}}$ is $(d - 2)-$dimensional (because $H_i$ and $H_j$ are not parallel). Denote by $\omega_{ij}$ the
unique unit norm vector orthogonal to both $H_{\{i,j\}}$ and $a_i$ (it is unique up to
reversal), and by $\Omega_{ij}$ be the 2-dimensional hyperplane through $\mu$ generated
by $a_i$ and $\omega_{ij}$.  By construction
$\Omega_{ij}$ is orthogonal to $H_{\{i,j\}}$ and is thus the orthogonal
complement to $H_{\{i,j\}}$ through $\mu$.  It follows that 
$\Omega_{ij} = \Omega_{ji}$ and by (\ref{important1}) 
that this hyperplane contains the orthogonal 
projection of $\mu$ onto $H_i$ (which is $y$) and the orthogonal 
projection of $\mu$ onto $H_j$ (which is $z$).  Further, 
$\Omega_{ij} \cap H_{\{i,j\}}$ consists of a single point $k_{\{i,j\}}$ satisfying
\begin{eqnarray}
\| \mu - k_{\{i,j\}} \|^2 & = & \| \mu - y \|^2 + \| y - k_{\{i,j\}} \|^2 \nonumber \\
& = & \| \mu - z \|^2 + \| z - k_{\{i,j\}} \|^2. \label{wedgie1}
\end{eqnarray}
Moreover
\begin{eqnarray}
&& y - k_{\{i,j\}} \, \mbox{is parallel to $\omega_{ij}$ and} \  z - k_{\{i,j\}} \, \mbox{is parallel to $\omega_{ji}$}, \nonumber\\
&& \| \mu - y \|^2 = \| \mu - z \|^2 = \rho, \ \ \mbox{and by (\ref{wedgie1})}, \nonumber \\
&& \| y - k_{\{i,j\}} \| = \| z - k_{\{i,j\}} \|, \ \ \mbox{and} \ \  \| \mu - y \| = \tan\phi \, \| y - k_{\{i,j\}} \|, \label{usestan}
\end{eqnarray}
where $2 \phi$ is the angle formed by $\omega_{ij}$ and $\omega_{ji}$. Let $h_{\{i,j\}}^g$ ($1 \le g \le d - 2$) be a basis for $\{ \, x \in \R^d \, : \, a_i^Tx = a_j^Tx = 0\}$.  Then $a_i$, together with $\omega_{ij}$ and the
$h_{\{i,j\}}^g$ form a basis for $\R^d$.  Let
\begin{itemize}
\item $O_i$ be the orthogonal projection of the origin onto $H_i$ -- hence $O_i$ is a multiple of $a_i$,
\item $N_i$ be the orthogonal projection of $O_i$ onto $H_{\{i,j\}}$.
\end{itemize}
We have
\begin{eqnarray}
y & = & O_i \, + \, (N_i - O_i) \, + \, (k_{\{i,j\}} - N_i) \, + \, (y - k_{\{i,j\}})\ 
\end{eqnarray}
\noindent and thus, since $N_i - O_i$ and $y - k_{\{i,j\}}$ are parallel to $\omega_{ij}$, and $k_{\{i,j\}} - N_i$ and $O_i$ are orthogonal to $\omega_{ij}$, 
\begin{eqnarray}
\omega_{ij}^T y & = & \omega_{ij}^T (N_i - O_i) + \omega_{ij}^T (y - k_{\{i,j\}}) \ = \ \omega_{ij}^T (N_i - O_i) + \| \omega_{ij} \| \| y - k_{\{i,j\}} \|,
\end{eqnarray}
or
\begin{eqnarray}
 \| y - k_{\{i,j\}} \| & = & \| \omega_{ij} \|^{-1} \omega_{ij}^T \, (y - N_i + O_i).
\end{eqnarray}
Consequently, by (\ref{usestan}) and (\ref{mudef})
\begin{eqnarray}
\hat \alpha(P^{i,j}, y) & = & \tan\phi \, \| y - k_{\{i,j\}} \| \nonumber \\
&& \nonumber \\
 & = & \tan\phi \, \| \omega_{ij} \|^{-1} \omega_{ij}^T \, (y - N_i + O_i), \label{hatalphadef}
\end{eqnarray}
which is affine as desired.  Note that for any 
$v \in H_{\{i,j\}}$ the lifting coefficient at $v$ as per (\ref{hatalphadef}) is
zero, 
 since $\omega_{ij}$ is orthogonal to both $a_i$ and $H_{\{i,j\}}$.
 Finally, 
since $\tan \phi$, $\omega_{ij}$, $N_i$ and $O_i$ are all polynomially computable the proof is now complete.   \QED \\

Now let $x^* \in \interior(P)$.  The problem
of finding the strongest possible lifted first-order inequality at $x^*$ chosen from among
those obtained by lifting from a point on the facet defined by the $i^{th}$ inequality can thus be written as follows:
\begin{eqnarray}
\mbox{{\bf SEP(i):}} \quad \min && - 2y^T x^* \, + \, \|y\|^2 \, - \, 2 \alpha (a_i^T x^* - b_i) \nonumber \\
s.t. && y \in P \nonumber \\
&& a_i^T y = b_i \nonumber \\
&& 0 \le \alpha \le p_{ij}^T y + q_{ij} \ \ \ \forall \ j \neq i. \nonumber
\end{eqnarray}
[Here, the last constraint is valid by the last part of the statement of
Lemma \ref{slice}.]
This is a linearly constrained, convex quadratic program with $d + 1$ variables and $2m - 1$ 
constraints.  

\begin{CO} Given $x^* \in \interior(P)$, in polynomial-time we can compute an
LFO inequality (\ref{lifted}) that attains the largest right-hand side
value at $x = x^*$.\end{CO}
\noindent {\em Proof.} We can attain the desired goal by solving SEP(i) for each choice of $1 \le i < m$. \QED

\subsection{Polynomial-time separation of LFOs when P is an ellipsoid and Q(x) is positive-definite quadratic}\label{ellipsoidal}
In this section we
will discuss a polynomial-time separation procedure for LFO inequalities in the case that $P$ is an ellipsoid with nonempty interior. As before we assume
without loss of generality that $Q(x) = \norm{x}^2$. Write
$$P \, = \, \{ x \, \in \R^d \, : \, x^T A x \, - \, 2 b^T x \, + \, c \, \le \, 0\}$$
for appropriate $A \succ 0$, $c$ and $b$.   To address the separation of LFO inequalities, consider a given a point $\bar x \in \text{int} (P)$. Using equation (\ref{LFOball}), the problem of finding an LFO inequality $q \ge \beta^T x + \beta_0$ whose right-hand side is maximized at $\bar x$ can be written as:
\begin{eqnarray}
q^{LFO} \ \doteq \ \min_{\mu, \rho} & & \norm{\mu}^2 - \rho - 2\bar x^T\mu \label{muopt} \\
  \mbox{s.t. } & & \{x \in \R^d : \norm{x - \mu}^2 \leq \rho\} \subseteq P \label{contained}\\
&& \mu \in \R^d, \quad \rho \ge 0. \label{boundsmuopt}
\end{eqnarray}
Since $P$ is compact, problem (\ref{muopt})-(\ref{boundsmuopt}) is well-posed, that is to say the objective is a ``min'' rather than simply an ``inf''. Denoting the maximum eigenvalue of $A$ by $\lambda_{max}$, we will
show in Corollary \ref{hatmuhatrho} given below that an optimal solution to this problem is given by
\begin{eqnarray}
\hat \mu & = & \lambda_{max}^{-1} b + (I - \lambda_{max}^{-1} A)\bar{x}, \nonumber \\
\hat \rho & = & \norm{\hat \mu - \bar x}^2 - \lambda_{max}^{-1}(\bar{x}^TA\bar{x} - 2b^T \bar{x} + c). \nonumber
\end{eqnarray}
A numerical application of this result to cardinality-constrained convex
quadratic programs will be provided in Section \ref{cardconst}.\\

For simplicity, in what follows we will
refer to problem (\ref{muopt})-(\ref{boundsmuopt}) as {\em the LFO separation problem at $\bar x$}, or the LFO separation problem for short.  Below we will provide a series of steps, culminating in Corollary \ref{hatmuhatrho}, that yield the description
of $\hat \mu$ and $\hat \rho$ given above.  First we present two technical 
results that will be used in the sequel.

\begin{LE} \label{getpi} Let $\tilde \mu \in \R^d$ and let $\tilde A$ a symmetric $d \times d$ matrix. The following two statements are equivalent:
\begin{eqnarray}
\mbox{{\bf (i)}} && \min_{x \in \R^d} \{ \, x^T \tilde A x - 2 \tilde \mu ^Tx \, \} > -\infty. \mbox{\hspace{.2in}}  \mbox{\hspace{1.5in}} \nonumber\\
\mbox{ {\bf (ii)}} && \mbox{$\tilde A \succeq 0$, and there exists $\pi \in \R^d$ satisfying } \ \tilde A \pi = \tilde \mu.\nonumber 
\end{eqnarray}
Furthermore, when (i) holds, a vector $\pi$ is an optimal solution to the minimization in (i) if and only if $\tilde A \pi = \tilde \mu$, in which case
\begin{eqnarray}
&& \min_{x \in \R^d} \{ \, x^T \tilde A x - 2 \tilde \mu ^Tx  \, \} \, = \, -\pi^T \tilde A \pi. \nonumber
\end{eqnarray}

\end{LE}
\noindent {\em Proof.}  Suppose (i) holds.   Then clearly $\tilde A \succeq 0$
and  for any $\delta \in \R^d$,
\begin{eqnarray}
 \quad \tilde \mu^T\delta = 0 \quad \mbox{whenever} \quad \tilde A \delta = 0. \nonumber
\end{eqnarray}
Farkas's Lemma then implies that there exists $\pi \in \R^d$ such that $\tilde A \pi = \tilde \mu$.  Thus (ii) holds;   the proof that (ii) implies (i) is similar
and will be omitted. Moreover, if (i) holds then (ii) does, and since the quadratic minimized in (i) is convex,  $x \in \R^d$ is an optimal solution to the minimization in (i) iff $ \tilde A x =  \tilde \mu$.  \QED \\

\begin{LE} \label{notquitethere} Let $\tilde v \in \R^d$
and let $\tilde A \succeq 0$ be a $d \times d$ matrix with maximum
eigenvalue $\tilde \lambda_{max} > 0$. Suppose $0 \le \theta \le \tilde \lambda_{max}^{-1}$. Then $x = \tilde v$ is an optimal solution to the problem
\begin{eqnarray}
\min_{x \in \R^d} \, \left\{ \, x^T(I -  \theta \tilde A)^T x  \ -2 \tilde v^T (I - \theta \tilde A) x  \, \right\}.
\end{eqnarray}
\end{LE}
\noindent {\em Proof.} Since $\theta \le \tilde \lambda_{max}^{-1}$ we have
$I - \theta \tilde A \succeq 0$.  The result now follows from Lemma \ref{getpi}
with $\tilde \mu = (I - \theta \tilde A) \tilde v$. \QED \\

\noindent We now return to the LFO separation problem at a given point $\bar x \in \interior(P)$; we will show that it is equivalent to:
\begin{eqnarray}
\hat q & \doteq & \min_{\mu, \theta} \,  \, \left[ -2\bar{x}^T\mu - \min_x\{x^T(I - \theta A)x - 2(\mu - \theta b)^Tx \} - \theta c \right] \label{slimobj} \\
&&    \text{s.t. }\quad \ \,  \mu \in \R^d, \quad 0 \le \theta \le \lambda_{max}^{-1}. \label{slimconst}
\end{eqnarray}
Our equivalence proof will proceed in several steps, given by Lemmas \ref{equiv1} and \ref{separationfromtheta}.  

\begin{LE} \label{equiv1} We have that $\hat q \le q^{LFO}$.  Further, suppose problem (\ref{slimobj})-(\ref{slimconst})
has an optimal solution $(\hat \mu, \hat \theta)$ such that $\hat \rho \ge 0$,
where
\begin{eqnarray}
&& \hat \rho \ \doteq \ \norm{\hat \mu}^2 + \min_x \{x^T(I - \hat \theta A)x - 2(\hat \mu - \hat \theta b)^Tx \} - \hat \theta c. \label{hatrhofromhat}
\end{eqnarray}
\noindent Then $\hat q = q^{LFO}$ and $(\hat \mu, \hat \rho)$ is optimal for the LFO separation problem.
\end{LE}
\noindent {\em Proof.} Since we assume $\interior(P) \neq \emptyset$, the S-Lemma (see \cite{slemma}, \cite{polter}, \cite{bennem}) implies that $(\mu, \rho)$ satisfies (\ref{contained}) if and only if there is some nonnegative real $\theta = \theta(\mu, \rho)$ such that
\[
\norm{x - \mu}^2 - \rho - \theta(x^TAx - 2b^Tx + c) \geq 0 \quad \forall x \in \mathbb{R}^d.
\]
This is equivalent to saying that there is $\theta \geq 0$ with
\begin{eqnarray*}
\min_x \{\norm{x - \mu}^2 - \rho - \theta(x^TAx - 2b^Tx + c)\} \ \geq \ 0,
\end{eqnarray*}
or equivalently
\begin{eqnarray}
\min_x \{x^T(I - \theta A)x - 2(\mu - \theta b)^Tx + \norm{\mu}^2 - \rho - \theta c\} \ \geq \ 0.  \label{lastmuopt}
\end{eqnarray}
We clearly must have $\theta \leq \lambda_{max}^{-1}$ for (\ref{lastmuopt}) to hold.  We can now write the LFO separation problem as:
\begin{eqnarray}
    \min_{\mu, \rho, \theta} &  & \norm{\mu}^2 - \rho - 2\bar x^T\mu \label{slemmaobj} \\
    \text{s.t. } & & \norm{\mu}^2 - \rho  + \min_x \{x^T(I - \theta A)x - 2(\mu -  \theta b)^Tx \} -  \theta c \ \geq \ 0 \label{slemmad}\\
&& \mu \in \R^d, \quad \rho \ge 0, \quad 0 \le \theta \le \lambda_{max}^{-1}. \label{slemmabounds}
\end{eqnarray}
\noindent We now use this formulation to argue that $\hat q \le q^{LFO}$.  Note that in any feasible
solution to (\ref{slemmaobj})-(\ref{slemmabounds})  the minimum in (\ref{slemmad}) must be finite and is therefore attained (as the quantity being minimized is
a quadratic). Hence, without loss of generality, (\ref{slemmad}) will hold with equality (or we could increase $\rho$); using this fact we can
eliminate $\norm{\mu}^2 - \rho$ from the objective. Thus $\hat q \le q^{LFO}$,
as desired.  Moreover if an optimal solution $(\hat \mu, \hat \theta)$ for problem (\ref{slemmaobj})-(\ref{slemmabounds}) is such that $\hat \rho$ defined as in (\ref{hatrhofromhat}) is
nonnegative, then clearly  $(\hat \mu, \hat \theta, \hat \rho)$ is 
feasible for (\ref{slemmaobj})-(\ref{slemmabounds}), with value at most $\hat q$. \hspace{.1in}\QED \\

\noindent Our next task is to further simplify problem (\ref{slimobj})-(\ref{slimconst}). Toward this goal we seek
an alternative description of the inner minimization in (\ref{slimobj}); 
in particular a characterization of those cases when it has finite value.  For $0 \le \theta \le \lambda_{max}^{-1}$ define
\begin{eqnarray}
&& F(\theta) \ \doteq \ \min_{\pi \in \R^d} \left\{ \, -2 \bar x^T \left[ \theta b + (I - \theta A) \pi \right]  
+ \pi^T(I -  \theta A)^T \pi \, \right\} \ + \ \theta c. \label{lastexpression}
\end{eqnarray}
\noindent We have:
\begin{LE} \label{separationfromtheta} (a) Problem (\ref{slimobj})-(\ref{slimconst}) can be equivalently rewritten as
\begin{eqnarray}
&& \min_{0 \le \theta \le \lambda_{max}^{-1}} F(\theta). \label{thetaobj}
\end{eqnarray}
(b) Further, suppose $\hat \theta$ is an optimal solution to (\ref{thetaobj}), 
and let $\hat \pi$ attain the minimum in (\ref{lastexpression}) when
$\theta = \hat \theta$.
If 
\begin{eqnarray}
\tilde \rho & \doteq & \norm{ \hat \mu - \hat \pi}^2 - \hat \theta ( \hat \pi^T A \hat \pi - 2 b^T \hat \pi + c) \ \ge \ 0 \label{haharho}
\end{eqnarray}
then
$(\hat \mu, \tilde \rho)$ is an optimal solution
for the LFO separation problem, where 
\begin{eqnarray}
&& \hat \mu = \hat \theta b + (I - \hat \theta A) \hat \pi. \label{hahamu} 
\end{eqnarray}
\end{LE}
\noindent{\em Proof.} We have $\hat q \le 0$, since a feasible
solution for the LFO separation problem is $\mu = \bar x$
and $\rho = 0$ with objective value $-\norm{\bar x}^2$.  It follows that
the inner minimum in (\ref{slimobj}) can be assumed to be finite. 
Thus, by Lemma \ref{getpi} we have that 
\begin{eqnarray}
&& \argmin_x \{x^T(I - \theta A)x - 2(\mu - \theta b)^Tx \} \, = \, \{ \pi \in \R^d \, : \,  (I - \theta A) \pi = \mu - \theta b \}. \label{argminnie}
\end{eqnarray}
Any $\pi$ in the set in (\ref{argminnie}) clearly satisfies 
\begin{eqnarray}
&& -\pi^T (I - \theta A) \pi - \theta c \, = \, \min_x \{x^T(I - \theta A)x - 2(\mu - \theta b)^Tx\} - \theta c. \nonumber
\end{eqnarray}
Consequently,
problem (\ref{slimobj})-(\ref{slimconst}) can be equivalently rewritten as
\begin{eqnarray}
&& \min_{\mu, \pi} \quad \ \  \ -2 \bar x^T \mu  \ + \ \pi^T(I - \theta A)\pi + \theta c  \label{fthetadefagain}\\
&& \quad \quad \quad \quad \  \mbox{s.t.} \ \  \mu - \theta b = (I - \theta A)\pi \label{restatemu}\\
&& \quad \quad \quad \quad \quad  \ \ \  \ \mu, \, \pi \in \R^d. \label{fthetabounds}
\end{eqnarray}
Substituting (\ref{restatemu}) into (\ref{fthetadefagain}) we obtain that
problem (\ref{thetaobj}) is equivalent to problem (\ref{fthetadefagain})-(\ref{fthetabounds}).  This proves part (a).  To prove (b), we have that
at $\theta = \hat \theta$ the optimum solution to (\ref{fthetadefagain})-(\ref{fthetabounds}) is $(\hat \mu, \hat \pi)$  
which satisfy equation (\ref{restatemu}), i.e. expression (\ref{hahamu}). 
Further, by Lemma \ref{equiv1}, if we can argue that
$$\hat \rho \ \doteq \ \norm{\hat \mu}^2 + \min_x \{x^T(I - \hat \theta A)x - 2(\hat \mu - \hat \theta b)^Tx \} - \hat \theta c \ \ge \ 0,$$
then $(\hat \mu, \hat \rho)$ is optimal for the LFO problem. But
using (\ref{argminnie}) and (\ref{hahamu}), $\hat \pi$ solves $\min_x \{x^T(I - \hat \theta A)x - 2(\hat \mu - \hat \theta b)^Tx\}$, and so
$$ \hat \rho = \norm{\hat \mu}^2 + \hat \pi ^T(I - \hat \theta A) \hat \pi - 2(\hat \mu - \hat \theta b)^T \hat \pi \,  - \, \hat \theta c \ = \tilde \rho$$
as defined in (\ref{haharho}), which is nonnegative by assumption. \QED \\

\noindent The following result characterizes $F(\theta)$.
\begin{LE} \label{notquitethere2} Let $0 \le \theta \le \lambda_{max}^{-1}$. 
Then $F(\theta)$ is obtained by choosing $\pi = \bar x$ in (\ref{lastexpression}).
Thus
\begin{eqnarray}
&& F(\theta) = -\norm{\bar x}^2 + \theta ( \bar x^T A \bar x - 2 b^T \bar x + c). \label{fthetasimple} 
\end{eqnarray}
\end{LE}
\noindent {\em Proof.}   The result follows from Lemma \ref{notquitethere},
with $\tilde v = \bar x$. \QED

\begin{CO}\label{hatmuhatrho} The optimizer for (\ref{thetaobj}) is $\hat \theta = \lambda_{max}^{-1}$, and an optimal solution to the LFO separation problem 
at $\bar x$ is
\begin{eqnarray}
\hat \mu & = & \lambda_{max}^{-1} b + (I - \lambda_{max}^{-1} A)\bar{x},  \label{hatmuoptimal} \\
\hat \rho & = & \norm{\hat \mu - \bar x}^2 - \lambda_{max}^{-1}(\bar{x}^TA\bar{x} - 2b^T \bar{x} + c). \label{hatrhooptimal}
\end{eqnarray}
\end{CO}
\noindent {\em Proof.}  By assumption $\bar x \in \interior(P)$.  
Hence the multiplier of $\theta$ in (\ref{fthetasimple}) is negative; 
consequently $F(\theta)$ is minimized at  $\hat \theta = \lambda_{max}^{-1}$.  
The result now follows from Lemmas \ref{notquitethere2} and \ref{separationfromtheta}, since the right-hand side of (\ref{hatrhooptimal}) is nonnegative because $\bar x \in \interior(P)$. \QED

\section{Tightening a general quadratic expression}\label{twoquad}
Consider a set of the form
\begin{eqnarray}
  && \Pi \ \doteq \ \{ \, (x, w, q) \in \R^n \times \R \times \R \ : \ q \geq x^THx + h^Tx, \ \ w \leq x^TAx \, \} \label{initPi}
\end{eqnarray}
where $H \succ 0$ and $A \succeq 0$ are $n \times n$ matrices.   With $P = \{(x, w) \in \R^n \times \R \, : \, x^TA x \leq w \},$ the set $\Pi$ is an example of our general set $S$ as in (\ref{Sdef}).  We will assume that $A$ has 
positive largest eigenvalue $\lambda_{max}$.
Here we show how the specialization of the LFO inequalities to this case
leads to a polynomial-time separable characterization of $\conv(\Pi)$. \\

\noindent As motivation for
this study, consider an optimization problem of the form
\begin{eqnarray}
&& \min \{ \, f(x) \, : \, x^T M x \, + \, h^T x \le  b_0, \,  x \, \in \, F \, \} \label{general1}
\end{eqnarray}
where $f \, : \, \R^n \rightarrow \R$, $M \in \R^n \times \R^n$ and $F \subseteq \R^n$, and $b_0 \in \R$.  We now apply a ``d.c.'' (difference between convex)
step (see \cite{dc}): we can always find matrices $H \succ 0$ and $A \succ 0$ such that
$x^T M x \ = \ x^T H x \, - \, x^T A x$ for all $x \in \R^n$.   Thus, (\ref{general1}) can be restated as
\begin{eqnarray}
&& \min \{ \, f(x) \, : \, q \ge x^T H x + h^T x, \ w \le x^T A x, \ q - w \le  b_0, \  x \, \in \, F \}, \nonumber
\end{eqnarray}
which can be relaxed to
\begin{eqnarray}
&& \min \{ \, f(x) \, : \, (x, w, q) \in \conv(\Pi), \ q - w \le  b_0, \  x \, \in \, F, \, q \in \R, \, w \in \R \, \}, \label{relax1}
\end{eqnarray}
where $\Pi$ is as in (\ref{initPi}).  A polynomial-time separation procedure
for $\conv(\Pi)$ can thus be used as a component in an algorithm for solving
the relaxation (\ref{relax1}).  \\

\noindent In the remainder of this section we address the 
separation problem for a set $\conv(\Pi)$ with $\Pi$ as in (\ref{initPi}).  We stress that we only require $A \succeq 0$, with positive maximum eigenvalue $\lambda_{max}$,
rather than $A \succ 0$ as in the above paragraph. In Section \ref{subsublfo} we first 
describe LFO inequalities as they pertain to the set $\Pi$ given above. 
In Section \ref{lfoparaboloid} we then
introduce another set of valid inequalities, the \textit{paraboloid} 
inequalities, which are valid for $\conv(\Pi)$.  We then show that, effectively, 
the
strongest paraboloid inequality at any given point is an LFO inequality. 
Finally, in Section \ref{nospectrum}  we show how to separate over paraboloid 
inequalities in polynomial time.

\subsubsection{LFO inequalities for\bmath{\conv(\Pi)}}\label{subsublfo}
Define $$P \doteq \{(x, w) \in \R^n \times \R_+ \, : \, x^TA x \leq w \},$$
When $A \succ 0$ (which is stronger than we require) this describes a paraboloid in $(x, w)$-space. With this notation,
$$ \Pi = \{ (x, w, q) \in \R^n \times \R \times \R \, : \, q \, \geq \, x^TH x + h^Tx, \ (x,w) \in \R^n \times \R - \, \interior(P) \},$$
and by Theorem \ref{domination} we have that $\conv( \Pi)$ is described by
linearization and LFO inequalities.  Since $\partial P$ is a differentiable
manifold homeomorphic to $\R^n$, it follows that $P$
is locally flat at every point $(\bar x, \bar x^T A \bar x) \in \R^n \times \R$, using as normal the unit vector in the
direction of
\begin{equation}
\begin{pmatrix} -2 A \bar x \\ 1 \end{pmatrix}. \nonumber
\end{equation}
Thus in the construction of the LFO inequality at $(\bar x, \bar x^T A \bar x)$ 
we consider inequalities of the form
\begin{equation} 
q \ \geq  \ \bar x^T H \bar x + h^T\bar x \ + \ \left(\begin{bmatrix}2 H \bar x + h\\ 0 \end{bmatrix} - \alpha \begin{bmatrix}2A\bar x\\ -1 \end{bmatrix}\right)^T \begin{bmatrix} x - \bar x\\ w - \bar x^T A \bar x \end{bmatrix},
\end{equation}
or in other words
\begin{eqnarray}
q & \ge & \bar x^T H \bar x + h^T\bar x + (2 H \bar x + h)^T(x - \bar x) \ + \  \alpha (- 2 \bar x^T A (x - \bar x)  \, + \, w - \bar x^T A \bar x ) \nonumber \\
& = &  -\bar x^T H \bar x \, + \, (2 H \bar x + h - 2 \alpha A \bar x)^T x + \alpha( \bar x^T A \bar x + w);\label{liftedcutorig}
\end{eqnarray}
when $\alpha = \hat \alpha$ we obtain the LFO inequality. 
As can be seen, (\ref{liftedcutorig}) strengthens 
the linearized inequality $q \ge \bar x^T H \bar x + h^T\bar x + (2 H \bar x + h)^T(x - \bar x)$ in the (excluded) region where $w > x^T A x$. \\

\noindent To simplify the discussion below we will assume, without loss of
generality, that $H = I$.  Thus the LFO inequality at $(\bar x, \bar x^T A \bar x)$ becomes
\begin{eqnarray}
q & \ge & (2\bar x + h - 2 \hat \alpha A \bar x )^T x + \hat \alpha w + \hat \alpha \bar x^T A \bar x - \norm{\bar x}^2. \label{liftedsimple0}
\end{eqnarray}

\subsubsection{Paraboloid inequalities}\label{lfoparaboloid}
In order to obtain a separation algorithm for inequalities 
(\ref{liftedsimple0}) we will first derive a geometrical characterization
of the set $P$, similar to that involving the ball inequalities developed in Section \ref{quadspecial}, equation (\ref{geometric}).  The intuitive reason that
paraboloid inequalities supersede ball inequalities is that the latter arose
in the context of our generic framework (\ref{Sdef}) with $Q(x)$ positive-definite.
Given $\mu \in \R^n$, $\nu \in \R_+$ and $\alpha > 0$, let 
\begin{eqnarray}
\Gamma^{\mu, \nu, \alpha} & \doteq & \{\, (x, w) \, \in \, \R^n \times \R \ : \  \norm{x - \mu}^2 + \alpha \nu  \ \le \ \alpha w \, \}, \label{truegamma1}
\end{eqnarray}
which defines a paraboloid. For technical reasons, for 
$\mu \in \R^n$, $\nu \in \R_+$ we also set
\begin{eqnarray}
\Gamma^{\mu, \nu, 0} & \doteq & \{\, (x, w) \, \in \, \R^n \times \R_+ \ : \  x = \mu, \ \nu \le w\} \label{truegamma2}
\end{eqnarray}
This is different from what the behavior of (\ref{truegamma1}) would achieve with $\alpha = 0$, which is the set $\{ (x, w)  \in  \R^n \times \R_+ \, : \, x = \mu \}$. On the other hand (\ref{truegamma2}) better reflects the limiting behavior
of (\ref{truegamma1}) at $x = \mu$ as $\alpha \rightarrow 0^+$.  Then it can be seen (proof omitted for brevity) that for any point $(x,w) \in \R^n \times \R$,
\begin{eqnarray}
(x,w) \in \R^n \times \R - \interior(P) \quad \mbox{iff} \quad (x,w) \in \R^n \times \R - \interior(\Gamma^{\mu, \nu, \alpha}), \ \ \mbox{for all $(\mu, \nu, \alpha)$ such that $\Gamma^{\mu, \nu, \alpha} \subseteq P$}. \nonumber
\end{eqnarray}
Using this characterization we have that for each triple $(\mu, \nu, \alpha) \in \R^n \times \R_+ \times \R_+$ with $\Gamma^{\mu, \nu, \alpha} \subseteq P$ the following inequality is valid for $\Pi$:
\begin{eqnarray}
\cP(\mu, \nu, \alpha):  && q \quad \geq \quad (2\mu + h)^Tx  + \alpha w \, - \, \norm{\mu}^2  \, - \, \alpha \nu, \nonumber
\end{eqnarray}
which cuts-off $\interior(\Gamma^{\mu, \nu, \alpha})$ as given by the following result. 
\begin{PR} Let  $(\mu, \nu, \alpha) \in \R^n \times \R_+ \times \R_+$,
and $(\bar x, \bar w) \in \R^n \times \R$. (a) Suppose $(\bar x, \bar w) \in \interior(\Gamma^{\mu, \nu, \alpha})$. If $\bar q \le \|\bar x\|^2 + h^T \bar x$ then $(\bar x, \bar w, \bar q )$ violates $\cP(\mu, \nu, \alpha)$.  (b)  
$(\bar x, \bar w, \|\bar x\|^2 + h^T \bar x)$ violates $\cP(\mu, \nu, \alpha)$ iff $(\bar x, \bar w) \in \interior(\Gamma^{\mu, \nu, \alpha})$.
\end{PR}
\noindent {\em Proof.}  (a) Since $(\bar x, \bar w) \in \interior(\Gamma^{\mu, \nu, \alpha})$, $-\alpha \bar w < -  \norm{\bar x - \mu}^2 - \alpha \nu$, and so
$$ 0 \, > \, \norm{\bar x}^2 - 2 \mu^T \bar x + \norm{\mu}^2 + \alpha(\nu - \bar w) \ \ge \ \bar q \ + \ \norm{\mu}^2 - (2\mu + h)^T \bar x  + \alpha( \nu - \bar w)$$
since $\bar q \le \|\bar x\|^2 + h^T \bar x$. (b) The proof of this fact is similar to that of (a) and will be omitted for brevity. \QED \\

\noindent We term $\cP(\mu, \nu, \alpha)$ a \textit{paraboloid} inequality, 
as an extension of
the ball inequalities (\ref{LFOball}).   We stress that paraboloid inequalities
are only defined for  triples $(\mu, \nu, \alpha) \in \R^n \times \R_+ \times \R_+$
such that $\Gamma^{\mu, \nu, \alpha} \subseteq P$.  In particular, we must have
$ \alpha \le \lambda_{max}^{-1}$,  where as stated before $\lambda_{max}$ is the largest eigenvalue of $A$.\\

\noindent For future reference, we state the following result which follows
directly from the definition of $\Pi$ and the paraboloid inequalities:

\begin{RE}\label{bindingpara} $\cP(\mu, \nu, \alpha)$ supports
$\Pi$ if and only if there exists $\bar x \in \R^n$ such that 
$$\norm{\bar x - \mu}^2 = \alpha(\bar x^T A \bar x - \nu),$$
in which case $\cP(\mu, \nu, \alpha)$ is valid for $\Pi$ and binding at $(\bar x, \bar x^T A \bar x, \norm{\bar x}^2 + h^T \bar x)$.
\end{RE}

\noindent In the next sequence of results we show that LFO and paraboloid
inequalities are essentially equivalent. We will first analyze those
paraboloid inequalities that are `strongest' at a given point, and then
consider paraboloid inequalities that are supporting for $\Pi$. The final proof will be provided in Theorem \ref{dom2} below.  \\

\noindent We will use the following notation:
$$U \doteq \{ \, (\mu, \nu, \alpha) \in \R^n \times \R_+ \times \R_+ \, : \, \Gamma^{\mu, \nu, \alpha} \subseteq P \, \}.$$

\begin{PR}\label{itsclosed} $U$  is closed.\end{PR}
\noindent {\em Proof.} Let $(\breve \mu, \breve \nu, \breve \alpha) \in \closure(U) - U$. Suppose first that $\breve \alpha > 0$.  
Since $\Gamma^{\breve \mu, \breve \nu, \breve \alpha} \not \subseteq P$, there
exists $(x, w)$ with 
$$ w \, < \, x^T A x \ \ \mbox{and}  \ \ \norm{x - \breve \mu}^2 + \breve \alpha \breve \nu \, \le \, \breve \alpha w.$$
It follows that we can find $\epsilon > 0$ and $\delta > 0$ such that
for any $(\mu, \nu)$ with $\norm{ \mu - \breve \mu} \le \delta$ and 
$|\nu - \breve \nu| \le \delta$ we nevertheless still have
$$ w + \epsilon \, < \, x^T A x \ \ \mbox{and}  \ \ \norm{x - \mu}^2 + \alpha \nu \, \le \, \alpha (w + \epsilon),$$
a contradiction.  The case $\breve \alpha = 0$ is similar and will be omitted.  \QED

\begin{PR} \label{strongestpara} Let $(\tilde x, \tilde w) \in \R^n \times \R$.  Then there is a paraboloid inequality 
whose right-hand side evaluated at $(\tilde x, \tilde w)$ is maximum  among all
paraboloid inequalities. \end{PR} 
{\em Proof.} The
right-hand side of $\cP(\mu, \nu, \alpha)$ at $(\tilde x, \tilde w)$ equals
$- \norm{ \tilde x - \mu}^2 + \norm{\tilde x}^2 + h^T \tilde x + \alpha \tilde w - \alpha \nu$; we want to maximize this expression subject to $(\mu, \nu, \alpha) \in U$.
Removing constants from the expression, this is equivalent to maximizing
\begin{eqnarray}
&& - \norm{ \tilde x - \mu}^2 + \alpha \tilde w - \alpha \nu \label{tomaximize}
\end{eqnarray}
subject to $(\mu, \nu, \alpha) \in U$.
A feasible choice for the maximization is $(\mu, \nu, \alpha) = (\tilde x, \tilde x^T A \tilde x, 0)$ (which is in $U$, by construction) for which (\ref{tomaximize}) attains value $0$.  Suppose 
\begin{eqnarray}
&& \kappa \ \doteq \ \sup \{\, - \norm{ \tilde x - \mu}^2 + \alpha \tilde w - \alpha \nu \ : \ 
(\mu, \nu, \alpha) \in U \, \} \ > \ 0. \label{itsclosed!}
\end{eqnarray}
We will argue that in this case when computing the supremum in (\ref{itsclosed!}) we can assume
that $(\mu, \nu, \alpha)$ can be constrained to lie in a bounded set; since
by Proposition \ref{itsclosed}, $U$ is closed, it will follow
that there the supremum is achieved, as desired.

So assume (\ref{itsclosed!}) holds. Now whenever $(\mu, \nu, \alpha) \in U$ we have $0 \le \alpha \le \lambda_{max}^{-1}$ and $\nu \ge 0$. Since $\kappa > 0$, nonnegativity of $\alpha$ implies that $\tilde w > 0$,  and since $\kappa > 0$
we have  $\nu \le \bar w$.  As a consequence
$\alpha \bar w - \alpha \nu$ is bounded (above and below) and so
$\norm{\mu}$ must be bounded as well.  This concludes the proof. \QED\\

\noindent Proposition \ref{strongestpara} motivates the following definition.

\begin{DE} Let $(\tilde x, \tilde w) \in \R^n \times \R$.  Then a paraboloid inequality 
whose right-hand side evaluated at $(\tilde x, \tilde w)$ is largest from among all
paraboloid inequalities is called a \textit{strongest} paraboloid inequality
at $(\tilde x, \tilde w)$. \end{DE}

\begin{LE}\label{tightpara} Suppose the inequality $\cP(\mu, \nu, \alpha)$ does not support $\Pi$.  Then $\alpha > 0$, and there exists $\epsilon > 0$ such that $\cP(\mu, \nu - \epsilon, \alpha)$ is a paraboloid inequality supporting $\Pi$.\end{LE}
\noindent {\em Proof.} If $\alpha = 0$ then by definition
$\Gamma^{\mu, \nu, \alpha} = \{(x, w) \in \R^n \times \R \, : \, x = \mu, \, w \geq \nu\}$ and trivially
$\cP(\mu, \nu, \alpha)$ holds as an equality at any point of the form 
$(\mu, w, \norm{\mu}^2 + h^T \mu)$.  Thus $\alpha > 0$. Since $\cP(\mu, \nu, \alpha)$ does not support $\Pi$, 
by Remark \ref{bindingpara} it follows that
\begin{eqnarray}
&& \norm{x - \mu}^2 - \alpha  x^T A x \ > \ -\alpha \nu, \ \ \forall \, x \in \R^n,  \label{meansnupositive}
\end{eqnarray}
from which we obtain $\nu > 0$.  The left-hand side of (\ref{meansnupositive}) is a quadratic which is lower bounded
by the quantity on the right-hand side. Thus, the minimum of the quadratic is attained
at some $x^0 \in \R^n$, and since at $x = \mu$ the quadratic takes nonpositive value, it does so at $x^0$ as well.  Therefore there exists $\epsilon > 0$ such that
$$ \norm{ x - \mu}^2 - \alpha  x^T A  x \ \ge \  - \alpha (\nu -\epsilon), \ \ \forall \, x \in \R^n,$$
with equality at $x^0$.   Since at $x^0$ the quadratic is nonpositive and $\alpha > 0$, we have $\nu - \epsilon \ge 0$.  It follows
that $P(\mu, \nu - \epsilon, \alpha)$ is a paraboloid inequality
 supporting $\Pi$. \hspace{.1in} \QED 
\begin{CO}\label{thestrongestpara}Let $(\bar x, \bar w) \in \R^n \times \R$. Then there is a paraboloid
inequality supporting $\Pi$ that is a strongest paraboloid
inequality at $(\bar x, \bar w)$. \end{CO} 
\noindent {\em Proof.} A strongest paraboloid inequality $\cP(\mu, \nu, \alpha)$
at $(\bar x, \bar w)$ exists by Proposition \ref{strongestpara}.  It
this inequality does not support $\Pi$, then by Proposition $\ref{tightpara}$
$\alpha > 0$ and $\cP(\mu, \nu - \epsilon, \alpha)$ supports $\Pi$ for
some $\epsilon > 0$.  But since $\alpha > 0$, $\cP(\mu, \nu - \epsilon, \alpha)$ is stronger at $(\bar x, \bar w)$ than $\cP(\mu, \nu, \alpha)$, a 
contradiction. \QED

\begin{PR}\label{trivialpara} Any LFO inequality is a paraboloid
inequality. \end{PR}
\noindent{\em Proof.}  Inequality (\ref{liftedsimple0}) is the same
as $\cP( \mu, \nu, \alpha)$ where $\alpha = \hat \alpha$, $\mu = \bar x - \alpha A \bar x$ and $\nu = - \bar x^T A \bar x + \alpha \norm{A \bar x}^2$. \QED \\

\noindent We can now prove the main result in this section.  We first remind the reader that the linearization inequality for the function $\norm{x}^2 + h^T x$ at any point $(\bar x, \bar w)$ is given by
$$ q \ \ge \ (2 \bar x + h)^T x - \norm{\bar x}^2.$$
\begin{THM}\label{dom2}
Let $(\bar x, \bar w, \bar q)$ be such that $(\bar x, \bar w) \in \interior(P)$
and $(\bar x, \bar w, \bar q)$ satisfies the linearization inequality at
$(\bar x, \bar w)$.
If $(\bar x, \bar w, \bar q)$ violates a paraboloid inequality, then a paraboloid
inequality that is maximally violated by $(\bar x, \bar w, \bar q)$ is
an LFO inequality. Conversely, any LFO inequality violated at 
$(\bar x, \bar w, \bar q)$ is a paraboloid inequality. \end{THM}  
\noindent {\em Proof.} Suppose that $(\bar x, \bar w, \bar q)$ violates some
paraboloid inequality $\cP(\mu, \nu, \alpha)$.  By Corollary \ref{tightpara} without loss of generality $\cP(\mu, \nu, \alpha)$ supports $\Pi$.  Thus by 
Theorem \ref{domination} there is an LFO inequality whose violation 
at $(\bar x, \bar w, \bar q)$ is at least as large as that of 
$\cP(\mu, \nu, \alpha)$.  But then by Proposition \ref{trivialpara} that
LFO inequality and $\cP(\mu, \nu, \alpha)$ are one and the same inequality.
The converse is similar and will be omitted. 
 \QED \\

The principal consequence of Theorem \ref{dom2} is that separation over the LFO
inequalities is equivalent to separation over the paraboloid inequalities. In
the next section we show how to do this in polynomial time.

\subsubsection{Polynomial-time separation of paraboloid inequalities}\label{nospectrum}  
In this section we provide a polynomial-time algorithm that, 
given a point $(\bar x , \bar w) \in \R^n \times \R$ with $\bar x^T A \bar x < \bar w$ computes a paraboloid inequality
$\cP( \mu^*, \nu^*, \alpha^*)$
that is strongest at $(\bar x,  \bar w)$.   By Proposition \ref{strongestpara}
such an inequality exists. We will refer to this task as {\em the paraboloid separation problem}.\\

\noindent In what follows we assume that the pair $(\bar x, \bar w)$ is given. Recall that inequality $\cP(\mu, \nu, \alpha)$ requires $ q \, \geq \, (2\mu + h)^Tx  + \alpha w  - \norm{\mu}^2   -  \alpha \nu$.  Therefore, the
paraboloid separation problem at $(\bar x, \bar w)$ can be stated as:
\begin{eqnarray}
    \min & & -(2\mu + h)^T\bar{x} - \alpha \bar w + \alpha \nu + \norm{\mu}^2 \label{inside1}\\
    \mbox{s.t. } && \Gamma^{\mu, \nu, \alpha} \quad \subseteq \quad \left\{(x, w) \in \R^n \times \R  :  x^TAx - w \le 0\right\} \label{contained2} \\
    && \mu \in \R^n, \, \nu \in \R_+, \, \alpha \ge 0. \label{bounds3}
\end{eqnarray}
We will show below (Corollary \ref{parmax}) that an optimal solution to this 
problem is:
\[
 \mu^* = (I - \lambda_{max}^{-1} A) \bar x, \quad \nu^* = -\lambda_{max} \norm{\mu}^2 \ + \ \bar x^T (\lambda_{max} I - A) \bar x, \quad \alpha^* = \lambda_{max}^{-1}.
\]
To obtain this result we will need a number of technical steps.

\begin{LE} \label{justabove} Let $(\hat \mu, \hat \nu, \hat \alpha)$ be an optimal solution to:
\begin{eqnarray}
    \min & &  - (2\mu + h)^T\bar{x} - \alpha \bar w +  \alpha \nu + \norm{\mu}^2  \label{poopiee1} \\
    \mbox{s.t. } & &\alpha \nu + \norm{\mu}^2 + \min_x \left\{ x^T(I - \alpha A)x - 2\mu^Tx\right\} \geq 0 \label{poopiee2} \\
    && \mu \in \R^n, \, \nu \in \R_+, \, \alpha \ge 0. \label{poopiee3}
\end{eqnarray}
Suppose $\hat \alpha > 0$. Then $(\hat \mu, \hat \nu, \hat \alpha)$  is an optimal solution
to problem (\ref{inside1})-(\ref{bounds3}).
\end{LE}
\noindent {\em Proof.} First note that the two objective functions are identical. We claim that for $\alpha > 0$ constraints (\ref{contained2}) 
and (\ref{poopiee2}) are identical.  This holds because (\ref{contained2}) holds if and only
if for all $x \in \R^n$
$$ \norm{x - \mu}^2 + \alpha \nu \ \ge \ \alpha x^T A x,$$
which is (\ref{poopiee2}).   Thus $(\hat \mu, \hat \nu, \hat \alpha)$ is
feasible for problem (\ref{inside1})-(\ref{bounds3}), and by the
preceding remarks, if it is not optimal for (\ref{inside1})-(\ref{bounds3}) 
then it is improved upon by a solution of the form $(\mu, \nu, 0)$.  But
any such vector is also feasible for (\ref{poopiee1})-(\ref{poopiee3}). 
Since we assumed  $(\hat \mu, \hat \nu, \hat \alpha)$ was optimal 
for (\ref{poopiee1})-(\ref{poopiee3}) the proof is complete. \QED \\

\noindent It seems intuitively clear that at optimality constraint (\ref{poopiee2}) will 
be binding. A formal proof is as follows:
\begin{PR} \label{wlogeq} Without loss of generality, at optimality 
for problem (\ref{poopiee1})-(\ref{poopiee3}), constraint (\ref{poopiee2})
will hold as an equality.\end{PR}
\noindent {\em Proof.} Let $(\hat \mu, \hat \nu, \hat \alpha)$ be optimal
for (\ref{poopiee1})-(\ref{poopiee3}).  If $\hat \alpha = 0$ the result
is clear (minimum in (\ref{poopiee2}) attained at $x = \hat \mu$). And
if $\hat \alpha > 0$, by Lemma \ref{justabove}, $\cP(\hat \mu, \hat \nu, \hat \alpha)$ is a strongest paraboloid inequality at $(\bar x, \bar w)$.  However,
if (\ref{poopiee2}) is not binding at $(\hat \mu, \hat \nu, \hat \alpha)$ 
then by Remark \ref{bindingpara}, $\cP(\hat \mu, \hat \nu, \hat \alpha)$ is not supporting.  This
is a contradiction by Lemma \ref{tightpara}. \QED \\

\noindent This result allows us to further simplify the paraboloid separation problem.
\begin{PR} \label{lastequiv0} Problem (\ref{poopiee1})-(\ref{poopiee3}) can be equivalently restated as:
\begin{eqnarray}
    \min & &  - 2\bar x^T(I - \alpha A) \pi  \ + \ \pi^T(I - \alpha A)\pi \ - \ \alpha \bar w \ - \  h ^T\bar{x} \label{lastobj} \\
\mbox{s.t.} &&  \pi \in \R^n, \, \alpha \ge 0. \label{lastRd}
\end{eqnarray}
Further, if $\pi^*, \alpha^*$ are optimal to this problem, then $(\mu^*, \nu^*, \alpha^*)$ are optimal for problem (\ref{poopiee1})-(\ref{poopiee3}), where
\begin{eqnarray}
\mu^* & = & (I - \alpha^*A) \pi^*, \nonumber \\
\alpha^* \nu^* & = & -\norm{\mu}^2 \, + \, \pi^{*T}(I - \alpha^*A)\pi^*. \nonumber
\end{eqnarray}
\end{PR}
\noindent {\em Proof.} By Proposition \ref{wlogeq} the paraboloid separation problem is equivalent to:
\begin{eqnarray}
    \min & &  - (2\mu + h)^T\bar{x} - \alpha \bar w \ - \ \min_x \left\{ x^T(I - \alpha A)x - 2\mu^Tx\right\} \nonumber \\
\mbox{s.t.}    && \mu \in \R^n. \nonumber
\end{eqnarray}
Since at optimality the inner minimum must be finite, we can now apply Lemma \ref{getpi}, with $\tilde A = I - \alpha A$ and $\tilde \mu = \mu$.  Lemma \ref{getpi} 
implies the already noticed fact that $\alpha \le \lambda_{max}^{-1}$, and it
implies that 
$$ \min_x \left\{ x^T(I - \alpha A)x - 2\mu^Tx\right\} = -\pi^T(I - \alpha A)\pi,$$ for some $\pi$ with $(I - \alpha A)\pi \ = \ \mu$.  Furthermore, we can
equivalently
rewrite the paraboloid separation problem as 
\begin{eqnarray}
    \min & &  - (2\mu + h)^T\bar{x} - \alpha \bar w \ + \ \pi^T(I - \alpha A)\pi \nonumber \\
\mbox{s.t. } &&  (I - \alpha A)\pi \ = \ \mu \nonumber \\
    && \pi \in \R^n, \ \mu \in \R^n,  \, \alpha \ge 0, \nonumber
\end{eqnarray}
\noindent which together with Proposition \ref{wlogeq} completes the proof. \QED

\begin{CO} \label{parmax} An optimal solution to the paraboloid separation problem is
\[
 \mu^* = (I - \lambda_{max}^{-1} A) \bar x, \quad \nu^* = -\lambda_{max} \norm{\mu}^2 \ + \ \bar x^T (\lambda_{max} I - A) \bar x, \quad \alpha^* = \lambda_{max}^{-1}.
\]
\end{CO}
\noindent {\em Proof.} Applying Lemma \ref{notquitethere}, we get that if 
$(\pi^*, \alpha^*)$ is an optimal solution to 
problem (\ref{lastobj})-(\ref{lastRd}),  we have $\pi^* = \bar x$. Thus the
optimal objective value of (\ref{lastobj})-(\ref{lastRd}) is
$$ -\bar x^T(I - \alpha^* A) \bar x - \alpha^* \bar w - h^T \bar x.$$
\noindent Since $ \bar x^T A \bar x - \bar w < 0$,  it follows that $\alpha^* = \lambda_{max}^{-1} > 0$ is the optimal
choice for $\alpha$, and using Propositions \ref{justabove} and \ref{lastequiv0} we obtain
the
 desired result. \QED

\section{Numerical experiments} \label{experiments}
In this section we present initial numerical experiments involving LFO cuts.  Our implementations are straightforward and in particular do not include any
cut management strategies.  Nevertheless the experiments are promising and
highlight interesting behavior of two problem classes. In Section \ref{cardconst} we consider cardinality-constrained convex quadratic programs; in 
Section \ref{disjunctcomparo} we compare LFO cuts with the disjunctive 
approach used in \cite{ceriasoares} in the context used
Section \ref{polyhedron}, i.e. $Q(x) = \norm{x}^2$ and $P$ is a polyhedron.

\subsection{Cardinality-constrained convex quadratic programs}\label{cardconst}
In this section we present preliminary experiments involving problems of the form
\begin{eqnarray}
&& \min \left\{ \, M(x) \ : \ x \in \Delta, \ \norm{x}_0 \le K \, \right\} \label{quadcard}
\end{eqnarray}
where $M(x)$ is a convex quadratic, $\Delta = \{ x \in \R^d_+ \, : \, \sum_j x_j = 1 \}$ is the unit simplex, $0 \le K \le d$,  and for $x \in \R^d$, $\norm{x}_0$ is the number of nonzero entries in $x$.  This problem class has been 
studied before, see e.g. \cite{danoq}, \cite{danoeig}, \cite{gunderoth} and
arises in several applications. When $M$ is positive-definite, $d$ is large and $K$ much smaller than $d$, problems of this type can be quite difficult.  
This is in particular the case if the solution to the relaxation
to (\ref{quadcard}) obtained by ignoring the cardinality constraint is contained in the (relative) interior of $\Delta$.   The goal of our experiments 
is to study the effect of using LFO cuts to obtain lower bounds on the 
value of problem (\ref{quadcard}).

Problem (\ref{quadcard}) can be formulated as a nonlinear mixed-integer 
program:
\begin{eqnarray}
\min && \ M(x) \label{mip1} \\
\mbox{s.t.} && \sum_j x_j = 1 \nonumber \\
&& x_j - y_j \le 1 \quad \mbox{and} \quad  y_j \in \{0,1\}, \quad \ 1 \le j \le d \nonumber \\
&& \sum_j y_j \le K, \quad x \ge 0. \label{miplast}
\end{eqnarray}

\noindent However this formulation can prove weak in difficult cases. A much stronger relaxation, the {\em perspective} 
relaxation, was used in \cite{gunderoth}; it is also related to the
disjunctive method in \cite{ceriasoares} (also see \cite{frangent}).  However, the perspective relaxation
can also prove computationally expensive; see \cite{danoeig}.

In our implementation of the LFO cuts we rely on the following result proved in \cite{danoeig} in a more general context. Here and below we denote $\cF \doteq \left\{ x \in \Delta \, : \, \norm{x}_0 \le K \, \right\}$, and 
for $\omega \in \R^d$ we write
$$ \rho(\omega) \ \doteq \ \frac{ (1 - \sum_{j \notin X}  \omega_j)^2}{K} \ \ + \ \ \sum_{j \in X} \omega^2_j,$$
where $X \subseteq \{1, \ldots, d\}$ is the set of indices of the $d - K$ smallest values $\omega_j$.

\begin{LE} \label{ballemma} Suppose $w \in \R^d$ satisfies $\sum_{j = 1}^d \omega_j \ = \ 1$. Then $ \min \{ \norm{y - \omega}^2 \, : \, y \in \cF \} = \rho(\omega).$ 
\end{LE}  

Using this result we can derive LFO cuts for problem (\ref{quadcard}).  Given
any point $\omega \in \Delta$, Lemma \ref{ballemma} guarantees that $\cF \subseteq \R^d - \interior(\cB(\omega, \sqrt{\rho(\omega)}))$
and the results in Section \ref{ellipsoidal} can be used to generate cuts.  We have implemented these ideas in the following straightforward 
cutting-plane procedure.\\

\noindent {\bf 0.} We initialize our formulation as $\min \{ q \, : \, q \ge M(x), \ x \in \Delta \}$. \\
\noindent {\bf 1.} Solve the current formulation, with solution $(\bar x, \bar q)$.  Let $P^{\bar x} \doteq \cB\left(\bar x, \sqrt{\rho(\bar x)}\right)$. \\
\noindent {\bf 2.} Compute an LFO cut that separates $(\bar x, \bar q)$ from the set $$\conv\left( \ \{ \, (x, q) \in \R^d \times R \, : \, q \ge M(x), \ x \in \R^d - \interior(P^{\bar x}) \, \right\} \ ) .$$
\noindent {\bf 3.} If no such cut is found, or if the violation of this cut by
$(\bar x, \bar q)$ is smaller than
a tolerance $\epsilon > 0$, or if the number of iterations exceeds a limit $T$, exit. Otherwise add the cut to the formulation
and return to step 1.\\

\noindent The cut in step 2 is obtained precisely as in Section \ref{ellipsoidal}, i.e. it is
a strongest LFO cut at $\bar x$.  In this implementation only one cut is
obtained from each given set $P^{\bar x}$.  The cutting procedure can be improved (see
the discussion following the numerical results).\\

\noindent In the results reported below, we
compare the strength of the lower bound obtained by our cutting-plane algorithm to two alternatives. The first one is the bound
obtained by running the mixed-integer
formulation (\ref{mip1})-(\ref{miplast}) using a commercial solver, with a very long time limit.  The second one is the bound obtained by 
an application of the S-Lemma to the cardinality-constrained problem (\ref{quadcard}).  This second approach was used in \cite{danoeig} and it can be summarized as follows: let
$$ x^* \ = \argmin \{ M(x) \, : \, e^T x = 1 \},$$
where $e = (1, \ldots, 1)^T.$
Then, as 
argued above
\begin{eqnarray}
&&  \min \{ M(x) \, : \, x \notin \interior( \cB( x^*, \sqrt{ \rho(x^*)} ) \} \label{Lstardef-first}
\end{eqnarray}
provides a valid lower bound to the value of problem (\ref{quadcard}).
Moreover, by optimality of $x^*$ we have 
$$ e^T \nabla M(x^*) = 0.$$  Consequently, by convexity, for any $y \in \R^d$
$$ M(x^* + y) \ \ge \ \lambda_{min} \norm{y}^2,$$
where $\lambda_{min}$ is the minimum eigenvalue of the Hessian of $M(x)$. Together
with (\ref{Lstardef-first}) 
we obtain that 
\begin{eqnarray}
&& L^* \doteq \ M(x^*) \, + \, \lambda_{min} \rho(x^*) \label{Lstardef}
\end{eqnarray}
is a lower bound to the value of problem (\ref{quadcard}).  This approach
can be viewed as an application of the S-Lemma.  As shown in 
\cite{danoeig}, $L^*$ improves not only on what the mixed-integer
programming formulation yields in practicable time, but, usually, on the value of the perspective relaxation as well.\\

\noindent In the experiments below,  
our cutting-plane algorithm using LFO 
inequalities was run with 
tolerance $\epsilon = 1.00 \times 10^{-3}$ and iteration limit $T = 10$. 
The problem instances considered in Tables \ref{cardtable-center} and \ref{cardtable-zero} were 
generated as follows.  In each case, the quadratic $M(x)$ is positive-definite
and separable.  Thus without loss of generality we can write
$M(x) = (x - x^0) ^T \Lambda (x - x^0)$ where $\Lambda = \diag\{\lambda_1, \ldots, \lambda_n\}$ and the $\lambda_i$ are positive, and $x^0 \in \R^d$.  Each $\lambda_i$ was chosen randomly, by drawing from the 
uniform distribution on the interval $ [1, 1 + \theta]$, 
where $\theta > 0$ is a fixed parameter.  Note that $x^0 = \argmin \{M(x) \, : \, x \in \R^d\}$. In our experiments we used $x^0 = 0$ and $x^0 = d^{-1} e$. \\

In the tables below, the
columns headed ``LFO-L'' and ``LFO-t'' describe the lower bound on
problem (\ref{quadcard}) and running time
produced the cutting-plane algorithm, respectively.  The mixed-integer programming formulation (\ref{mip1})-(\ref{miplast}) was run until either a limit of 1000 CPU seconds was reached or one million branch-and-cut nodes were enumerated (whichever came first); columns headed ``MIP-L'', ``MIP-t'' and
``MIP nodes'' indicate, respectively, the resulting lower bound, (wall-clock) running time on sixteen threads, and number of 
branch-and-cut nodes. ``MIP-U'' provides the {\em upper} bound on problem
(\ref{quadcard}) obtained by the mixed-integer programming approach.  Finally,
the column headed ``S-L'' is the lower bound provided by the S-Lemma approach as
in (\ref{Lstardef}).  On the problem instances with $d > 100$ the objective functions were scaled up
by a factor of $1000$.

All computations (here and in the next section) were performed on an 8-core i7 computer, with 48 GB of physical memory.  The cutting-plane algorithm used Gurobi 5.50 \cite{gurobi} for step 1.  To run
the mixed-integer programs we used both Gurobi 5.50 and CPLEX 12.2 \cite{cplex} (and report the better of the two). 

\begin{table}[h]
\centering
\caption{{\bf \emph{Cardinality-constrained problems with\bmath{x^0 = d^{-1}e}}}}
\begin{tabular}{|r|r|r||r|r|r|r||r|r||r|}
\hline \hline
{\bf d} &{\bf K}& \bmath{\theta} & {\bf LFO-L}& {\bf S-L} & {\bf MIP-L} & {\bf MIP-U}& {\bf LFO-t}& {\bf MIP-t} &  {\bf MIP}\\
 & & & & & &  & (sec)& (sec) & {\bf nodes}\\
\hline \hline
100 & 20 & 2.00 &  0.0411 & 0.0412  & 0.0005 &  0.0587   & 0.127  & 227  & 1011704   \\
100 & 50 & 5.00 &  0.0108 &  0.0108 & 0.0006 & 0.0314    & 0.102  & 222  & 1004975   \\
100 & 20 & 10.00 & 0.0465   &0.0465   & 0.0009 & 0.1284    & 0.120  & 288  &  1008679  \\
1000 & 100 & 10.00 &  9.1009 & 9.1010  & 0.0010 & 18.2534    & 0.883  & 1012  & 246063   \\
1000 & 100 & 100.00 &  10.0109  & 10.0125  & 0.0048  & 87.8492    & 0.848  & 1004  & 208633   \\
1000 & 70 & 20.00 & 13.5842  & 13.5844 & 0.0011 & 32.0741    & 0.879  & 1000  & 176152   \\
2000 & 100 & 40.00 & 9.5178 & 9.5178  & 0.0003 & 26.8787    & 3.014  & 1086  & 34699   \\
2000 & 90 & 50.00 & 10.6348  &  10.6358  & 0.0003 & 32.2729    & 2.563  & 1019  & 14298   \\
2000 & 80 & 50.00 & 12.0266  & 12.0280    & 0.0003 & 33.8795    & 3.186  & 1015  & 152638   \\
\hline
\end{tabular}
\label{cardtable-center}
\end{table}

Table \ref{cardtable-center} reports on results using $x^0 = d^{-1} e$.  Problems of this type are especially 
hard for the mixed-integer programming formulation, which is unsuccessful at moving the lower bound significantly
away from zero.  This also holds for the smaller problem instances, even though over a million branch-and-cut nodes
are enumerated. The LFO-based approach quickly (within ten iterations) attains a bound that is essentially equal to that provided by the S-Lemma,
and several orders of magnitude larger than the mixed-integer lower bound, and
thereafter tails off, sharply.  The bound proved by the LFO-based approach, in many cases, does not 
completely close the gap relative to the best upper bound obtained by the mixed-integer
programming formulation in the provided time/node limit.

To further explore this behavior consider Table \ref{cardtable-zero} which displays results in cases where
$x^0 = 0$.  When this is the case one can prove that the optimal value of problem (\ref{quadcard}) is obtained
as follows: where $I \subseteq \{1, \ldots, d\}$ is the set of indices corresponding to the $K$ smallest $\lambda_i$,
the optimal value of problem (\ref{quadcard}) equals $ \left( \sum_{j \in I} \lambda_j^{-1} \right)^{-1}.$
Table \ref{cardtable-zero} displays this value in the column headed ``OPT''.  We can see that the mixed-integer 
solver obtains this value as an upper bound (but does not prove so) nearly all the time (in one case, highlighted
with an asterisk, round-off error by the solver resulted in a better-than-optimum upper bound).  In this case,
again, we see that the lower bound obtained by using the mixed-integer formulation is nearly always greatly improved
by the LFO-driven lower bound (which as before essentially ties the S-Lemma lower bound).  Generally, the LFO
lower bound reduces the duality gap by at least $50\%$.

\begin{table}[h]
\centering
\caption{{\bf \emph{Cardinality-constrained problems with\bmath{x^0 = 0}}}}
\begin{tabular}{|r|r|r||r|r|r|r|r||r|r||r|}
\hline \hline
{\bf d} &{\bf K}& \bmath{\theta} & {\bf LFO-L}& {\bf S-L} & {\bf MIP-L} & {\bf MIP-U}& {\bf OPT} & {\bf LFO-t}& {\bf MIP-t} &  {\bf MIP}\\
 & & & & & &  & & (sec)& (sec) & {\bf nodes}\\
\hline \hline
100 & 50 & 0.50 & 0.0204  & 0.0204 & 0.0127  & 0.0224 & 0.0224 & 0.11  & 260 & 1016899 \\
100 & 20 & 0.50 & 0.0489  & 0.0490 & 0.0129 & 0.0524 & 0.0491  & 0.10 & 258 & 1009229\\
100 & 10 & 0.50 &0.0984  & 0.0984 & 0.0103 & 0.1026 & 0.1026  & 0.09 & 250 & 1011881\\
100 & 50 & 1.00 & 0.0213  & 0.0213 & 0.0157  & 0.0246 & 0.0214  & 0.09  & 309  & 1022007  \\
100 & 20 & 1.00 & 0.0484  & 0.0485 & 0.0151  & 0.0548 & 0.0487  & 0.09  & 293  & 1013290  \\
100 & 10 & 1.00 & 0.0972  &  0.0972& 0.0158  & 0.1053 & 0.1053  & 0.12  & 242  &1004596  \\
100 & 50 & 4.00 & 0.0288  &0.0288   &  0.0254 & 0.0362 & 0.0362 & 0.13  & 330   & 1001651  \\
1000 & 100 & 10.00 & 10.0295  & 10.0322   & 4.2512  &  14.7429 & 14.7429   & 1.01  & 1013   & 171452  \\
1000 & 90 & 10.00 & 10.9764  & 10.9763   &  4.2736 &  15.8986 &  15.8986 & 0.86  &  1015  &  170395 \\
1000 & 80 & 10.00 & 12.1727  & 12.1760   &  4.2541 &  17.3035 & 17.3012   & 0.979  &  1005  &  169706 \\
1000 & 70 & 20.00 & 14.7284  & 14.7285   & 6.7351  &  23.2463$^*$ &   23.2464  & 0.638  & 1005   & 170092 \\
2000 & 100 & 40.00 & 11.0015  &  11.0038  & 5.5414  & 18.9675 & 18.9675  & 2.738  & 1011  & 111417  \\
2000 & 90 & 50.00 & 12.5227  &  12.5229  &  6.5581 & 22.1117 & 22.1117 & 2.168  & 1025  & 103016  \\
2000 & 80 & 50.00 &13.5913   & 13.5914    & 6.5594  & 23.6386 &  23.6386  & 2.267  & 1056  &   107836  \\
\hline
\end{tabular}
\label{cardtable-zero}
\end{table}

The above cutting-plane scheme would likely be improved in a number of ways. Principal among
these is the concept of {\em sampling} the infeasible region so as to generate strong cutting-planes in advance of the
formal algorithm, possibly used in preprocessing form.  A number of such sampling techniques (related to the so-called Sandwich algorithm) are described
in \cite{bao}.  The dampened method in Step 3 of Example \ref{ex2} can also be viewed as an example of this idea.  
Preprocessing by sampling so as to generate good cuts in advance of the formal algorithm is usually a very effective
idea, especially in the context of first-order algorithms used to approximate a very nonlinear function.  

\noindent In forthcoming work we plan to address the following  more substantial  enhancements to the preliminary work described here:
\begin{itemize}
\item[(i)] The approach given by steps {\bf 0 - 4} above
relies on Lemma \ref{ballemma} to exclude a ball $\cB$ from the feasible region, with the current iterate $\bar x$ at
its center.   Typically, the optimizer $z$ of the quadratic $M(x)$ on the boundary of this ball will be such that
$z - \bar x$ is parallel to an eigenvector corresponding to the minimum eigenvalue for $M(x)$.  However, the ball
$\cB$, computed as per Lemma \ref{ballemma}, has as radius the minimum distance from $\bar x$ to the feasible region --
and this minimum distance is attained by some point $y \in \partial \cB$ such that $y - \bar x$ is, often, \textit{far from} parallel  with any eigenvector arising from a small eigenvalue for $M(x)$.  In fact  $y - \bar x$ may more 
likely be aligned with eigenvectors corresponding to much larger eigenvalues.  
This phenomenon suggests the use a polyhedral, rather than a spherical, relaxation.  Rather than excluding a 
ball from the feasible region, we would
exclude a polyhedron which ``pushes'' into the corners of the unit simplex.  Ideally, 
the excluded polyhedron would have larger diameter along the
eigenspace corresponding to the larger eigenvalues.
\item[(ii)] Couple the cutting-plane scheme with a branching approach, 
parameterized by a value $\theta > 0$, as follows. Let $v$ be an eigenvector
of the quadratic part of $M(x)$ corresponding to the maximum eigenvalue.  If $\hat x$ is the solution to the 
current relaxation, we would branch by considering those points $x$ such that $v^T(x - \hat x) > \theta$, 
those such that $v^T(x - \hat x) < -\theta$, and those such that $ -\theta \le v^T(x - \hat x) \le \theta$ -- and on the
latter set use a relaxation that enforces a minimum distance to the feasible region which would presumably be 
larger than the minimum distance from $\hat x$ to the feasible region.
\item [(iii)] Reformulate the problem so as to use a coordinate system contained in the $(d-1)$-dimensional hyperplane $\{ x \in \R^d \, : \, e^T x = 0\}$.
As shown in \cite{danoeig}, an appropriate representation of the quadratic
$M(x)$
valid in this restricted space (the so-called projected quadratic) results in more
effective bounds, because the minimum eigenvalue usually increases.  This point is
related to (i), above.
\end{itemize}

\subsection{A comparison with the disjunctive method}\label{disjunctcomparo}
In this section we consider the setup in section \ref{polyhedron}.  We are given a polyhedron $P = \{ x \in \R^d \, : \, Ax \leq b \}$ containing the origin in its interior and we are interested in the set $S = \{(x,q) \in \R^d \times \R \, : \, q \ge \norm{x}^2, \ x \in \R^{d} - \interior(P) \}$.

The purpose of the experiments in this section is to compare the performance of the separation algorithm for LFO
inequalities given in Section \ref{polyhedron} to the performance of the disjunctive method, both 
as generators of cutting-plane used to separate from $\conv(S)$. Here we remind the reader that a set $S$ of the
form considered here would not arise as ``the'' problem being solved. Rather, it would be a \textit{relaxation} of
a problem of interest (as was the case with the cardinality-constrained problem considered above).  Thus, our
algorithm for separating LFO inequalities on the one hand, and the disjunctive formulation (plus SOCP duality)
on the other hand, would constitute competing methods for separating from $\conv(S)$, and here we investigate
their performance in this context.

In the case of either separation routine we run a cutting-plane algorithm similar to the one given in
the Section above:\\

\noindent {\bf 0.} We initialize our formulation as $\min \{ q \, : \, q \ge M(x) \}$. \\
\noindent {\bf 1.} Solve the current formulation, with solution $(\bar x, \bar q)$.  \\
\noindent {\bf 2.} Compute a cut that separates $(\bar x, \bar q)$ from the set $$\conv\left( \ \{ \, (x, q) \in \R^d \times R \, : \, q \ge M(x), \ x \in \R^d - \interior(P) \, \right\} \ ) .$$
\noindent {\bf 3.} If no such cut is found, or if the violation of this cut by
$(\bar x, \bar q)$ is smaller than
a tolerance $\epsilon > 0$, or if the number of iterations exceeds a limit $T$, exit. Otherwise add the cut to the formulation
and return to step 1.\\

This algorithm solves the problem $N_2 \doteq \min \{ \norm{x}^2 \, : \, x \notin \interior(P) \}$.  Of course the value of
this problem is known, however here we are concerned with the number of iterations needed by the two algorithms and with any other
interesting performance attributes that may arise.\\

In our experiments, the systems $Ax \leq b$ were generated as follows. The entries of the coefficient vector $a_i$ were set to random uniform values between $-1$ and $1$, and then each was set to $0$ with probability $0.5$. The vector was rejected if it was a positive multiple of any of the previous vectors $a_1, \hdots, a_{i-1}$. $a_i$ was then normalized to have unit norm and the entries were rounded to $3$ digits. Next, the value $\bar b_i$ was calculated as
\[
\bar b_i = \max\{a_i^Tx \ \vert \ a_j^Tx \leq b_j, \ j = 1, \hdots, i-1\}.
\]
If $\bar b_i$ was finite, we set $b_i$ to a value randomly distributed between $0.5\bar b_i$ and $0.95\bar b_i$. Otherwise $b_i$ was set to $1 + \Gamma$, where $\Gamma$ was a generated randomly from a gamma distribution with shape $\sqrt{n}$ and scale $0.5\sqrt{n}$. In either case $b_i$ was then rounded to $3$ digits.

In Table \ref{outpoly-1},  the columns labeled $d$, $m$, and $val$ give the dimension, number of rows in $A$, and true problem value $N_2$. $Lo_l$, $Time_l$, and $Cuts_l$ give the best lower bound, time taken (in seconds), and number of cuts generated by the lifting method, with similar information for the disjunctive method given by 
$Lo_d$, $Time_d$, and $Cuts_d$, respectively. Each method, using
the template provided by steps {\bf 0-4} above, was allowed a total of $500$ cuts, and was only allowed to generate a new cut if fewer than 600 seconds had passed since the initial setup. An asterisk next to a time in the $time_d$ column indicates that that instance was stopped because the solver was unable to find the dual variables needed to generate a cut. These tests were terminated when the relative gap between the lower bound and the true value was less than a tolerance of $10^{-5}$.  

\begin{table}[h]
\centering
\caption{{\bf \emph{Comparison with disjunctive method}}}
\begin{tabular}{|r|r||r|r|r||r|r||r|r|}
\hline
{\bf d} &{\bf m}& \bmath{N_2} & \bmath{Lo_{l}}&\bmath{Lo_{d}}&\bmath{Time_{l}}&\bmath{Time_{d}}&\bmath{Cuts_{l}}&\bmath{Cuts_{d}}\\
\hline \hline
10 & 50 & 5.191 & 5.190 & 2.121 & 0.2   & 74.7& 51 & 500\\
\hline
20 & 100 & 14.537 & 14.536 & 0.366 & 0.1  & 88.1& 15  & 500\\
\hline
20 & 300 & 17.831 & 17.830 & 0.771 & 0.3  & 311.4& 15  & 500\\
\hline
50 & 200 & 79.888 & 79.887 & 0.154 & 0.3  & 381.5 & 20 & 500\\
\hline
75 & 250 & 343.897 & 336.353 & 36.872 & 126.7  & 317.1$^*$ & 500  & 95\\
\hline
100 & 300 & 324.486 & 324.485 & 16.139 & 0.6  & 126.7$^*$ & 14 & 11\\
\hline
200 & 400 & 2207.060 & 2207.038 & 0.000 & 8.7  & 91.0$^*$  & 92 & 1\\
\hline
300 & 500 & 4583.748 & 4583.733 & 0.000 & 2.4  & 155.3$^*$& 20  & 1\\
\hline
800 & 1200 & 38142.592 & 38142.243 & 0.000 & 32.1 & 1879.1$^*$ & 42  & 1\\
\hline
1000 & 2000 & 61726.150 & 61725.542 & 0.000 & 227.8   & 3304.3$^*$& 134 & 1\\
\hline
\end{tabular}
\label{outpoly-1}
\end{table}

We can see from this table that the disjunctive method tends to fail as problem size becomes large.  Part of the reason is that the SOCPs to be solved simply prove too difficult.  Note that the disjunctive method already seems to encounter
numerical difficulties on medium size problems;  even on the smallest instances, the lower bounds
on $N_2$ obtained by the disjunctive method can be poor.

Table \ref{outpoly-2} presents results on various variants of the cutting-plane algorithm. 
In addition to the LFO cuts, the ``Basic" version does not use the constraint $q \ge \norm{x}^2$ and instead uses linearization cuts to approximate this constraint. Its initial formulation is 
$\min \{ q \, : \, q \ge 0 \}$.  We also consider enhancements to the Basic version, using three heuristics to help make faster progress:
\begin{enumerate}
\item Before starting, the linearization cut was added at each unit vector $e_i$ as well as $-e_i$.  
\item Before starting, the LFO inequality at the point closest to the origin on each facet was added, if possible.
\item The constraint $Ax \leq b$ was added in the relaxation.
\end{enumerate}
The first two heuristics are versions of ``sampling'' as described at the end of the last section.  The ``Full" method includes the conic constraint $q \ge \norm{x}^2$ and does not use linearization cuts.
These tests were terminated if, between subsequent iterations, the objective value $z$ and all entries of the solution $x$ were within a tolerance of $10^{-3}$ of the previous values.  
In Table \ref{outpoly-2}, columns labeled $d$, $m$, and $N_2$ are as in the previous table. In the next three sets of columns (``Basic", ``Heuristics", and ``Full"), $q_{lo}$ gives the best proven lower bound, $lin$ gives the number of linearization cuts added, $lifted$ gives the number of lifted cuts added, and $t$ is the time spent, in seconds. In this test, each method was limited to $30$ minutes to add cuts. During each iteration, the relaxation was solved and both the linearization cut and the lifted cut were added, if possible. A maximum of $10,000$ iterations were reached.

We can see from Table \ref{outpoly-2} that on all instances but the last one, the ``Full'' version is the clear winner.
The use of the conic constraint helps guide the algorithm toward points where cutting (using LFO inequalities) is
most effective.  At the same time, having a single conic constraint helps control computational cost and numerical
instabilities.

\begin{landscape}
\begin{table}[h]
\centering
\caption{{\bf \emph{Comparison of cutting-plane strategies}}}
\begin{tabular}{| r|r|r || r|r|r|r || r|r|r|r || r|r|r |}
 \hline
\multicolumn{3}{|c||}{} & \multicolumn{4}{c||}{{\bf Basic}} & \multicolumn{4}{|c||}{{\bf Heuristics}} & \multicolumn{3}{c|}{{\bf Full}}\\
 \hline
\textbf{$d$} &\textbf{$m$} & \textbf{$N_2$} & \textbf{$q_{lo}$}& \textbf{${lin}$} & \textbf{${lifted}$}& \textbf{$t$}& \textbf{$q_{lo}$}& \textbf{${lin}$} & \textbf{${lifted}$}& \textbf{$t$}& \textbf{$q_{lo}$}& \textbf{${lifted}$}&\textbf{$t$}\\
\hline \hline
10 & 50 & 5.19 & 5.19 & 286 & 120 & 0.8 & 5.19 & 280 & 212 & 1.1 & 5.19 & 133 & 1.1 \\
\hline
20 & 100 & 14.54 & 14.54 & 2104 & 824 & 19.5 & 14.54 & 890 & 744 & 8.6 & 14.54 & 33 & 0.3 \\
\hline
20 & 300 & 17.83 & 17.83 & 1903 & 577 & 28.7 & 17.83 & 2831 & 13 & 34.0 & 17.83 & 26 & 0.6 \\
\hline
50 & 200 & 79.89 & 0.00 & 10000 & 1 & 688.0 & 79.89 & 10300 & 38 & 646.2 & 79.89 & 28 & 0.6 \\
\hline
75 & 250 & 343.90 & 0.00 & 10000 & 1 & 768.6 & 343.87 & 8198 & 7841 & 1803.9 & 339.46 & 1370 & 1802.3 \\
\hline
100 & 300 & 324.49 & 0.00 & 10000 & 1 & 1127.2 & 320.57 & 6780 & 6361 & 1804.1 & 324.49 & 25 & 1.1 \\
\hline
200 & 400 & 2207.06 & 0.00 & 269 & 1 & 1807.1 & 2207.06 & 6168 & 151 & 1803.3 & 2207.03 & 87 & 9.6 \\
\hline
300 & 500 & 4583.75 & 0.00 & 301 & 1 & 13.3 & 597.60 & 2938 & 2003 & 1803.3 & 4583.74 & 24 & 3.7 \\
\hline
800 & 1200 & 38142.59 & 0.00 & 302 & 1 & 66.7 & 569.48 & 2940 & 651 & 1811.4 & 38142.54 & 55 & 53.9 \\
\hline
1000 & 2000 & 61726.15 & 0.00 & 246 & 1 & 108.7 & 61726.15 & 4001 & 893 & 2330.6 & 57689.25 & 18 & 40.0 \\
\hline
\end{tabular}
\label{outpoly-2}
\end{table}

\end{landscape}

\subsubsection{Example}
Consider the bilinear form
$$ f(x) \ \doteq \ 2( x_1 x_2 \, + \, x_1 x_3 \, + x_2 x_3 )$$
over the unit cube $[0,1]^3$.  Writing, for $1\le i < j \le 3$, $f_{ij} = x_i x_j$, the McCormick relaxation for $f_{ij}$ amounts to:
$$ f_{ij} \, \ge \, x_i + x_j - 1, \ \ \ f_{ij} \le \min\{ x_i , x_j \}.$$
At $\bar x = (1/2 , 1/2, 1/2)^T$, the lower bound on $f(\bar x)$ produced by the McCormick relaxation is zero (for more complex examples see \cite{jeff1}).  We show next how
our procedures may be used to generate a formulation that proves a positive
lower bound on $f(\bar x)$.  We stress that what we have here is an ad hoc
construction -- we plan to return to this topic in a future work.

\noindent We have $f(x) = U(x) - L(x)$, where
\begin{eqnarray}
U(x) & \doteq & (x_1 + x_2)^2 + (x_1 + x_3)^2 + (x_2 + x_3)^2, \nonumber \\
L(x) & \doteq & 2(x_1^2 + x_2^2 + x_3^2). \\
\end{eqnarray}
Now we apply the techniques from Section \ref{twoquad}.  We have
$U(x) = x^T H x$ and $L(x) = x^T A x$, where
\[
H = \begin{bmatrix} 2&1&1\\1&2&1\\1&1&2 \end{bmatrix}, \quad\quad A = \begin{bmatrix} 2&0&0\\0&2&0\\0&0&2\end{bmatrix}.
\]
The Cholesky decomposition of $H$ is 
\[
H = LL^T = \begin{bmatrix} \sqrt{2}&0&0\\1/\sqrt{2}&\sqrt{3/2}&0\\1/\sqrt{2}&1/\sqrt{6}&\sqrt{4/3} \end{bmatrix} \begin{bmatrix} \sqrt{2}&1/\sqrt{2}&1/\sqrt{2}\\0&\sqrt{3/2}&1/\sqrt{6}\\0&0&\sqrt{4/3} \end{bmatrix}.
\]
Let $V\Lambda V^T$ be the eigendecomposition of $L^{-1}AL^{-T} = 2 L^{-1}L^{-T} = 2 (L^TL)^{-1}$:
\[
V = \begin{bmatrix}\sqrt{3}/6&1/2&2/\sqrt{6}\\1/6&-3/\sqrt{2}& \sqrt{2}/3\\ -2\sqrt{2}/3&0&1/3 \end{bmatrix}, \quad\quad \Lambda = \begin{bmatrix} 2&0&0\\0&2&0\\0&0&1/2\end{bmatrix}.
\]
The transformation we use is $p = V^TL^Tx$, or $x = L^{-T}Vp$. Note:
\[
L^{-T}V = \begin{bmatrix}1/\sqrt{2}&-1/\sqrt{6}&-\sqrt{3}/6\\0&2/\sqrt{6}&-3/\sqrt{6}\\0&0&\sqrt{3}/2\end{bmatrix}\begin{bmatrix}\sqrt{3}/6&1/2&2/\sqrt{6}\\1/6&-3/\sqrt{2}& \sqrt{2}/3\\ -2\sqrt{2}/3&0&1/3 \end{bmatrix} = \begin{bmatrix}1/\sqrt{6}& \sqrt{2}/2& \sqrt{3}/6\\1/\sqrt{6}&-\sqrt{2}/2&\sqrt{3}/6\\-2/\sqrt{6}&0&\sqrt{3}/6 \end{bmatrix}.
\]
Thus, we have
\[
x \in [0, 1]^3 \Leftrightarrow \begin{bmatrix}-L^{-T}V\\ L^{-T}V\end{bmatrix}p \leq \begin{bmatrix}0\\0\\0\\1\\1\\1 \end{bmatrix}.
\]
Let $\breve H$ be the image of $[0,1]^3$ under the mapping. It can be seen
that for any $x$ we have 
$$p_3(x) = \frac{2}{\sqrt{3}} \, (x_1 + x_2 + x_3)$$
and thus our point of interest, $\bar x$, is mapped to $\bar p = (0, 0, \sqrt{3})^T$.

\noindent Further, in $p$-space, $f(x)$ is represented as
$$ F(p) \ \doteq \ (p_1^2 + p_2^2 + p_3^2) \  - \ (2p_1^2 + 2p_2^2 + \frac{1}{2} p_3^2).$$

\noindent Consider the ``paraboloid cut''
\begin{eqnarray}
p_1^2 + p_2^2 + (p_3 - 2\alpha\sqrt{3})^2 + \epsilon & \geq & 2p_1^2 + 2p^2_2 + \frac{1}{2}p_3^2. \label{par1}
\end{eqnarray}
For $\alpha = \epsilon = 1/10$, a calculation shows that (\ref{par1}) is 
valid for all $p \in \breve H$ with $p_3 \ge \sqrt{3}$ (or, informally,
it is valid for all $x \in [0,1]^3$ with $\sum_i x_i \ge 3/2.$)  In the
region of validity, we therefore have
$$ F(p) \ \ge \ 4\alpha\sqrt{3}p_3 - 12\alpha^2 - \epsilon \ = \ \frac{2}{5} \sqrt{3} p_3 - \frac{11}{50}.$$
In other words, for $x \in [0,1]^3$ with $\sum_i x_i \ge 3/2$,
$$ f(x) \ \ge \ \frac{4}{5} (x_1 + x_2 + x_3) - \frac{11}{50}.$$
Consider now the paraboloid cut (\ref{par1}) with $\alpha = 1/2$ and $\epsilon = 3/2$.  A calculation shows that in that case (\ref{par1}) is 
valid for all $p \in \breve H$ with $p_3 \le \sqrt{3}$.  Where it is valid we
get 
$$ F(p) \ \ge \ 4\alpha\sqrt{3}p_3 - 12\alpha^2 - \epsilon \ = \ 2\sqrt{3} p_3 -  \frac{9}{2},$$
and thus, for $x \in [0,1]^3$ with $\sum_i x_i \le 3/2$,
$$ f(x) \ge 4(x_1 + x_2 + x_3) - \frac{9}{2}.$$
\noindent We now have a disjunction between two polyhedra:
\begin{eqnarray}
 \Theta &\doteq & \left\{ (x, f) \ : \ x \in [0,1]^3, \ \sum_j x_j \ge 3/2, \ f \ge \frac{4}{5} (x_1 + x_2 + x_3) - \frac{11}{50} \right\}, \ \ \mbox{and} \nonumber\\
\Pi &\doteq & \left\{ (x, f) \ : \ x \in [0,1]^3, \ \sum_j x_j \le 3/2, \ f \ge 4(x_1 + x_2 + x_3) - \frac{9}{2} \right\}. \nonumber
\end{eqnarray}
Thus, solving the linear program
\begin{eqnarray}
&& \min \, f \nonumber \\
\mbox{s.t.} && (x,f) \in \conv( \Theta \cup \Pi) \nonumber \\
  && x = \bar x \nonumber
\end{eqnarray}
yields a valid lower bound on $f(\bar x)$.  The value of this LP is $.43599$.
Using LP duality, one also obtains the valid cut
$$ f(x) \ \ge \ 1.162667(x_1 + x_2 + x_3) \ - \ 1.308001,$$
which proves the lower bound $f(\bar x) \ge .43599$.  

As the example makes 
clear, numerical precision issues should be an integral component of an implementation of the ideas presented in this paper.  We plan to take up these
points in a future work.\\

\noindent{\bf Acknowledgments.} We would like to thank Samuel Burer, Adam Letchford and G\'{a}bor Pataki for helpful suggestions.

\tiny  Tue.Dec.17.095932.2013

\begin{thebibliography}{99}
\bibitem{aldon} {\sc F. Alizadeh and D. Goldfarb}, Second-order cone
programming, {\em Mathematical Programming} {\bf 95} (2001), 3 -- 51.

\bibitem{alkha} {\sc F. Al-Khayyal and J. Falk}, Jointly constrained biconvex programming, {\em Math. Oper. Res.} {\bf 8} (1983), 273-­286.

\bibitem{andetal} {\sc K. Andersen, Q. Louveaux, R. Weismantel and L. Wolsey}, Cutting planes from two rows of a simplex tableau, in : {\em Integer Programming and Combinatorial Optimization, 
Lecture Notes in Computer Science} {\bf 4513} (2007), pp 1 -­ 15.


\bibitem{kurtalone} {\sc K.M. Anstreicher}, Semidefinite programming versus the reformulation-linearization technique for nonconvex quadratically constrained quadratic programming, {\em J. Global Optimization} {\bf 43} (2009), 471--484.


\bibitem{anstreicherburer} {\sc K.M. Anstreicher and S. Burer}, Computable representations for convex hulls of low-dimensional quadratic forms, {\em Mathematical Programming B} {\bf 124} (2010), 33--43.

\bibitem{dc} {\sc K.M. Anstreicher and S. Burer}, D.C. versus copositive bounds for standard QP, {\em Journal of Global Optimization} {\bf 33} (2005), 299 -­ 312.


\bibitem{alper} {\sc A. Atamt\"{u}rk and V. Narayanan}, Lifting for conic mixed-integer programming, {\em Mathematical Programming} {\bf 126} (2011), 351--363. 


\bibitem{bao} {\sc X. Bao, N.V. Sahinidis, and M. Tawarmalani}, Multiterm polyhedral relaxations for nonconvex,
quadratically constrained quadratic programs,  {\em Optimization Methods and Software} {\bf 24} (2009), 485--504.

\bibitem{bal71} {\sc E. Balas}, Intersection cuts - a new type of cutting planes for integer programming, {\em Operations Research} {\bf 19} (1971), 19 -­ 39.


\bibitem{bal75} {\sc E\@. Balas}, Disjunctive programs: cutting planes from
logical conditions, in O.L. Mangasarian et al., eds., {\em Nonlinear 
Programming} {\bf 2} (1975), Academic Press, 279 -- 312.

\bibitem{bal79} {\sc E\@. Balas}, Disjunctive programming. {\em Annals
of Discrete Mathematics} {\bf 5} (1979), 3 -- 51.

\bibitem{balcercor93} {\sc E\@. Balas, S\@. Ceria and G\@. Cornu\'{e}jols},
A lift-and-project cutting plane algorithm for mixed 0-1 programs,
{\em Mathematical Programming} {\bf 58} (1993), 295 -- 324.

\bibitem{bennem} {\sc A. Ben-Tal and A. Nemirovsky}, {\em Lectures on Modern Convex Optimization: Analysis, Algorithms, and Engineering Applications} (2001),
MPS-SIAM Series on Optimization,  SIAM, Philadelphia, PA.

\bibitem{couenne} {\sc P. Belotti, J. Lee, L. Liberti, F. Margot and A. Wachter}, Branching and bounds tightening techniques for non-convex
   MINLP, {\em Optim. Methods Softw.} {\bf 24} (2009), 597-­634.


\bibitem{belmillnam}{\sc P. Belotti, A.J. Miller and M. Namazifar}, Valid inequalities and convex hulls for multilinear functions, {\em Electronic Notes in Discrete Mathematics} {\bf 36} (2010), 805 -­ 812.

\bibitem{danoq}{\sc D. Bienstock}, Computational study of a family of mixed-integer quadratic programming problems, {\it Mathematical Programming} {\bf 74} (1996), 121-- 140.


\bibitem{danoeig} {\sc D. Bienstock}, Eigenvalue techniques for proving bounds for convex objective, nonconvex programs, in: 
{\em Integer Programming and Combinatorial Optimization, 
Lecture Notes in Computer Science} {\bf 6080} (2010), pp 29--42. 

\bibitem{boydvan} {\sc S. Boyd and L. Vandenberghe}, {\em Convex Optimization},
Cambridge University Press (2004).

\bibitem{bucalo} {\sc C. Buchheim, A. Caprara and A. Lodi}, An effective branch-and-bound algorithm for convex quadratic programming, 
{\em Integer Programming and Combinatorial Optimization, 
Lecture Notes in Computer Science} {\bf 6080} (2010), 285 -- 298.


\bibitem{burerletchford} {\sc S. Burer and A. N. Letchford}, On non-convex quadratic programming with box constraints, {\em SIAM Journal on Optimization} {\bf 20} (2009), 1073--1089.

\bibitem{burlet} {\sc S. Burer and A. N. Letchford}, Non-convex mixed-integer nonlinear programming: a survey, {\em Optimization Online}, February 2012.



\bibitem{lee}{\sc S. Cafieri, J. Lee, and L. Liberti}, On convex relaxations of quadrilinear terms, {\em Journal of Global Optimization} {\bf 47} (2010), 661-–685.

\bibitem{ceriasoares} {\sc S\@. Ceria and J. Soares}, Convex programming for disjunctive convex optimization, {\em Mathematical Programming} {\bf 86} (1999), 595 -­ 614.

\bibitem{cornmarg} {\sc G. Cornu\'{e}jols and F. Margot}, On the facets of mixed integer programs with two integer variables and two constraints, {\em Mathematical Programming} {\bf 120} (2009), 429­-456.


\bibitem{cplex} {\sc IBM ILOG CPLEX Optimizer}. http://www-01.ibm.com/software/integration/optimization/cplex-optimizer/.


\bibitem{ismael1} {\sc I.R.\ de Farias JR., E.L.\ Johnson and G.L.\ Nemhauser}, Facets of the complementarity knapsack polytope, {\em Mathematics of Operations Research} {\bf 27}  (2002), 210 --226.

\bibitem{lifttilt} {\sc D.G. Espinoza, R. Fukasawa and M. Goycoolea}, Lifting, tilting and fractional programming revisited, {\em Oper. Res. Lett.} {\bf 38} (2010), 559 -- 563.

\bibitem{frangent} {\sc A. Frangioni and C. Gentile}, Perspective cuts for a class of convex 0-1 mixed
integer programs, {\em Mathematical Programming} {\bf 106} (2006), 225 -­ 236.



\bibitem{dongarud} {\sc D. Goldfarb and G. Iyengar}, Robust portfolio selection problems, {\em Mathematics of Operations Research} {\bf 28}  (2002), 1 -- 38.


\bibitem{golub1} {\sc G.H. Golub}, Some modified matrix eigenvalue problems, {\em SIAM Review} {\bf 15} (1973), 318 -- 334.

\bibitem{gunderoth} {\sc O. G\"{u}nl\"{u}k and J. T. Linderoth}, Perspective reformulations of mixed integer nonlinear programs with indicator variables, {\em Mathematical Programming} {\bf 124} (2010), 183 --205.

\bibitem{gurobi} {\sc Gurobi Optimizer}. http://www.gurobi.com/.

\bibitem{ismael2} {\sc A.B.\ Keha, I.R.\ de Farias JR. and G.L.\ Nemhauser}, A branch-and-cut algorithm without binary variables for nonconvex piecewise linear optimization,'' {\em Operations Research} {\bf 54} (2006), 847 -- 858.

\bibitem{jeff1} {\sc M. Kilinc, J. Linderoth and J. Luedtke}, Effective separation of disjunctive cuts for convex mixed integer nonlinear programs, {\em Optimization Online} (2010).

\bibitem{kojtunc1} {\sc M. Kojima and L. Tuncel}, Discretization and 
localization in successive convex relaxation methods for nonconvex quadratic optimization problems, {\em Mathematical Programming} {\bf 89} (2000), 79 -- 111.

\bibitem{lacilex} {\sc L\@. Lov\'{a}sz and A\@. Schrijver},
Cones of matrices and set-functions and 0-1 optimization,
{\em SIAM J. on Optimization} {\bf 1} (1991), 166-190.

\bibitem{lnl} {\sc J. Luedtke, M. Namazifar and J. Linderoth}, Some results on the strength of relaxations of multilinear functions, Optimization Online, August 2010.

\bibitem{alphastrong} {\sc J. Mak\'{o}, K. Nikodem and Z. P\'{a}les}, 
On strong $(\alpha, \F)$-convexity, {\em Mathematical Inequal. Appl.} {\bf 2} (2012),
289 -- 299.

\bibitem{mccormick} {\sc G.P. McCormick}, Computability of global solutions to factorable nonconvex programs: Part I {$-$} Convex underestimating problems.
{\em Mathematical Program.} {\bf 10} (1976), 147 -- 175.

\bibitem{mitchelletal} {\sc J.E. Mitchell, J.-S. Pang and B. Yu},
Convex quadratic relaxations of nonconvex quadratically constrained quadratic
programs, {\em Optimization Methods and Software}, (published online) 2012.

\bibitem{alex2012} {\sc A\@. Michalka}, Ph.D. Dissertation, Columbia University (in preparation).

\bibitem{nemwols} {\sc G.L. Nemhauser and L.A. Wolsey}, {\it Integer and
Combinatorial Optimization}, Wiley, New York (1988).


\bibitem{polter} {\sc I. P\'{o}lik and T. Terlaky}, A survey of the S-lemma, 
{\em SIAM Review} {\bf 49} (2007), 371 -- 418.

\bibitem{qbmarg} {\sc A. Qualizza, P. Belotti, F. Margot},  Linear programming relaxations of quadratically constrained quadratic programs, manuscript, 2011.

\bibitem{richtaw} {\sc J.-P. Richard and M. Tawarmalani}, Lifting inequalities: a framework for generating strong cuts for nonlinear programs, 
{\em Mathematical Programming} {\bf 121} (2010), 61 -- 104.


\bibitem{rikun} {\sc A.D. Rikun}, A convex envelope formula for multilinear functions. {\em Journal of Global Optimization} {\bf 10} (1997), 425 -­ 437.

\bibitem{baron} {\sc N.V. Sahinidis}, BARON: A general purpose global optimization software package,  {\em J. Global Opt.} {\bf 8} (996), 201 -­ 205.


\bibitem{sax1} {\sc A. Saxena, P. Bonami and J. Lee}, Convex relaxations of non-convex mixed integer quadratically constrained programs: Extended formulations. {\em Mathematical Programming B} {\bf 124} (2010), 383 -- 411.

\bibitem{sax2} {\sc A. Saxena, P. Bonami and J. Lee}, Convex relaxations of non-convex mixed integer quadratically constrained programs: Projected formulations. To appear, {\em Mathematical Programming}.

\bibitem{SA1}{\sc H. D. Sherali and W. P. Adams}, {\em A Reformulation-Linearization Technique for Solving Discrete and Continuous Nonconvex Problems}, Kluwer, Dordrecht (1998).

\bibitem{SA2} {\sc S\@. Sherali and W\@. Adams}, Reformulation-linearization technique (RLT) for semi-infinite and convex programs under mixed 0-1 and general discrete restrictions, {\em Discrete Applied Mathematics} {\bf 157} (2009), 1319--1333.

\bibitem{sheraliadams} {\sc S\@. Sherali and W\@. Adams},
A hierarchy of relaxations between the continuous and convex
hull representations for zero-one programming problems,
{\em SIAM J. on Discrete Mathematics} {\bf 3} (1990), 411-430.


\bibitem{stubbsmehrotra} {\sc R.A. Stubbs and S. Mehrotra}, 
A branch-and-cut method for 0-1 mixed convex programming, {\em Mathematical Programming} {\bf 86} (1999), 515 -­ 532.

\bibitem{talsah3} {\sc M. Tawarmalani and N.V. Sahinidis}, 
Convex extensions and envelopes of lower semi-continuous functions, {\em Mathematical Programming} {\bf 93} (2002), 247 -- 263.

\bibitem{talsah1} {\sc M. Tawarmalani and N.V. Sahinidis}, 
Global optimization of mixed-integer nonlinear programs: a theoretical and computational study, {\em Mathematical Programming} {\bf 99} (2004), 563 -- 591.

\bibitem{talsah2} {\sc M. Tawarmalani and N.V. Sahinidis}, 
A polyhedral branch-and-cut approach to global optimization, {\em Mathematical Programming} {\bf 103} (2005), 225 -­ 249.


\bibitem{slemma} {\sc V. A. Yakubovich}, S-procedure in nonlinear control theory, Vestnik Leningrad University, {\bf 1} (1971), 62--77. 

\end{thebibliography}
\end{document}